\documentstyle[12pt]{article}
\def\Bbb R{{\rm \bf R}}
\def\proclaim#1{\vskip2mm{\bf #1}\em}
\def\endproclaim{\em \vskip2mm}
\def\tag#1{\eqno(#1)}
\def\gathered{\begin{array}{c}}
\def\endgathered{\end{array}}
\def\text{\mbox}

\begin{document}

\title {The enclosure method for inverse obstacle scattering problems
with dynamical data over a finite time interval: III.  Sound-soft
obstacle and bistatic data}
\author{Masaru IKEHATA\footnote{
Department of Mathematics,
Graduate School of Engineering,
Gunma University, Kiryu 376-8515, JAPAN}}
\maketitle

\begin{abstract}
This paper is concerned with an inverse obstacle problem which
employs the dynamical scattering data of acoustic wave over a {\it
finite time interval}. The unknown obstacle is assumed to be
sound-soft one. The governing equation of the wave is given by the
classical wave equation. The wave is generated by the initial data
localized outside the obstacle and observed over a finite time
interval at a place which is not necessary the same as the support of the initial
data. The observed data are the so-called {\it bistatic data}. 
In this paper, an {\it enclosure method} which employs the bistatic
data and is based on two main analytical formulae, is developed.
The first one enables us to extract the maximum {\it spheroid}
with focal points at the center of the support of the initial data
and that of the observation points whose exterior encloses the
unknown obstacle of general shape. The second one, under some
technical assumption for the obstacle including convexity as an
example, indicates the deviation of the geometry of the boundary
of the obstacle and the maximum spheroid at the contact points.
Several implications of those two formulae are also given. In
particular, a constructive proof of a uniqueness of a {\it
spherical} obstacle using the bistatic data is given.

\noindent
AMS: 35R30, 35L05, 35J05

\end{abstract}


\section{Introduction}

In this paper, we consider an inverse obstacle scattering problem
for a {\it sound-soft} obstacle with {\it dynamical data} over a
{\it finite time interval}. The governing equation of the wave is
the classical wave equation.  The wave as the solution is
generated by the initial data whose support is localized at the
outside of the obstacle and observed over a finite time interval
on a different position from the support of the initial data. The
observed data are the so-called {\it bistatic data}.
This is a simple mathematical model of the data collection process using an acoustic wave/electromagnetic wave
such as, bistatic active {\it sonar}, {\it radar}, etc.
See, e.g., \cite{C} for the bistatic active sonar.
The aim of
this paper is to develop an enclosure method which employs the
bistatic data.

Let us describe a mathematical formulation of the problem.
Let $D$ be a nonempty bounded open subset of $\Bbb R^3$ with $C^2$-boundary
such that $\Bbb R^3\setminus\overline D$ is connected.
Let $0<T<\infty$.
Let $f\in\,L^2(\Bbb R^3)$ satisfy $\text{supp}\,f\cap\overline D=\emptyset$.
Let $u=u_f(x,t)$ denote the weak solution of the following initial boundary value problem for the classical wave equation:
$$\begin{array}{c}
\displaystyle
\partial_t^2u-\triangle u=0\,\,\text{in}\,(\Bbb R^3\setminus\overline D)\times\,]0,\,T[,\\
\\
\displaystyle
u(x,0)=0\,\,\text{in}\,\Bbb R^3\setminus\overline D,\\
\\
\displaystyle
\partial_tu(x,0)=f(x)\,\,\text{in}\,\Bbb R^3\setminus\overline D,\\
\\
\displaystyle
u=0\,\,\text{on}\,\partial D\times\,]0,\,T[.
\end{array}
\tag {1.1}
$$
Here $\nu$ denotes the unit {\it outward normal} to $D$ on $\partial D$.
The boundary condition for $u$ in (1.1) means that $D$ is
a {\it sound-soft} obstacle.  In this paper, $T$ is always fixed.  Thus, for our purpose
the weak solution over the bounded interval $]0,\,T[$ is appropriate.
Since the notion of the weak solution
for the wave equation is well established, we do not repeat the description here.
Instead see \cite{DL} for the notion
and also its use in \cite{IE0, IE3} for inverse obstacle scattering problems with dynamical data
over a finite time interval.


In this paper, we consider the following problem.

{\bf\noindent Inverse Problem.}
Let $B$ and $B'$ be two {\it known} open balls centered at $p\in\Bbb R^3$ and $p'\in\Bbb R^3$ with radius $\eta$
and $\eta'$, respectively such that $\overline B\cap\overline D=\emptyset$
and $\overline B'\cap \overline D=\emptyset$.
Let $\chi_B$ denote the characteristic function of $B$ and set $f=\chi_B$.
Assume that $D$ is {\it unknown}.
Extract information about the location and shape of $D$ from the data $u_f(x,t)$ given at all
$x\in B'$ and $t\in\,]0,\,T[$.

As far as the author knows,
there is no result to this problem for general configulation of $B$ and $B'$.
This is the problem raised in \cite{IE3} as an open problem related to the enclosure method itself.
In particular, the problem contains the case when $\overline B\cap\overline{B'}=\emptyset$ which 
corresponds to the case when the emitter and receiver are placed
on {\it different positions} at a {\it finite distance} from
the obstacle.  Strictly speaking, we should call the data in this case the bistatic data,
however, we include also the case $\overline B\cap\overline{B'}\not=\emptyset$.

In this paper, we develop an enclosure method
with bistatic data.  
In short, the enclosure method aims at extracting a domain that encloses an unknown discontinuity,
such as cavities, cracks, inclusions or obstacles.
The idea of the enclosure method goes back to \cite{I1}, in which the original enclosure method was
developed by considering an inverse boundary value problem
governed by the Laplace equation.
In \cite{I4}, an idea for the application of the enclosure method to the dynamical data
coming from the heat or wave equations has been introduced.
Now we have many applications of this enclosure method to inverse
boundary value problems governed by the heat equations in
\cite{IK1, IK2, IE2}, visco-elastic system of equations \cite{IH} and inverse obstacle scattering problems
governed by the wave equations in \cite{IE0, IE3, IE4}.

We establish two main analytical formulae.  
The first one enables us to extract the maximum {\it spheroid}
with focal points at the center of the support of the initial data
and that of the observation points whose exterior encloses the
unknown obstacle of general shape.  The appearence of the exterior
of a spheroid as an enclosing domain is {\it new} since previous enclosing domains
are a half plane/space, sphere or its exterior, or cone.
The formula shows us an effect of the bistatic data on the obtained information.
See Theorem 1.1 below.
The second one, under some
technical assumption for the obstacle including convexity as an
example, indicates the deviation of the geometry of the boundary
of the obstacle and the maximum spheroid at the contact points.
This is also new.
See Theorem 1.3 below.
And also we present several implications of those two formulae. In
particular, we give a constructive proof of a uniqueness of a {\it
spherical} obstacle using the bistatic data.

\subsection{Extracting the first reflection distance and its implication}

In this paper, given an arbitrary $h\in L^2(\Bbb R^3)$,
we denote by $v_h$ the unique weak solution $v\in H^1(\Bbb R^3)$ of
$$\displaystyle
(\triangle-\tau^2)v+h(x)=0\,\,\text{in}\,\Bbb R^3.
\tag {1.2}
$$
$v_h$ has the expression
$$\displaystyle
v_h(x)=v_h(x,\tau)=\frac{1}{4\pi}
\int_{\Bbb R^3}\frac{e^{-\tau\vert x-y\vert}}{\vert x-y\vert}h(y)dy.
\tag {1.3}
$$

Define
$$\displaystyle
w_f(x)=w_f(x,\tau)
=\int_0^Te^{-\tau t}u_f(x,t)dt,\,\,x\in\Bbb R^3\setminus\overline D,\,\,\tau>0.
$$
$w=w_f$ satisfies
$$\begin{array}{c}
\displaystyle
(\triangle-\tau^2)w+f(x)=e^{-\tau T}F_f(x,\tau)\,\,\text{in}\,\Bbb R^3\setminus\overline D,
\\
\\
\displaystyle
w=0\,\,\text{on}\,\partial D,
\end{array}
\tag {1.4}
$$
where
$$\displaystyle
F_f(x,\tau)=\partial_tu_f(x,T)+\tau u_f(x,T),\,\,x\in\Bbb R^3\setminus\overline D.
$$
Since $F_f(x,\tau)$ is unknown, it seems that the existence of such term in (1.4)
hides the information about an unknown obstacle. 
However, the use of the enclosure method presented below does not make 
it a problem at all and enables us to
extract the information about the obstacle provided $T$ is sufficiently large and fixed.

Let $\chi_{B'}$ denote the characteristic function of $B'$
and set $g=\chi_{B'}$.

The results of this paper are concerned with the asymptotic
behaviour of the {\it indicator function}:
$$\displaystyle
\tau\longmapsto
\int_{\Bbb R^3\setminus\overline D}(fv_g-w_fg)dx
=\int_Bv_gdx-\int_{B'}w_fdx
$$

For the description of the results we prepare some notation.

Define
$$\displaystyle
\phi(x;y,y')
=\vert y-x\vert+\vert x-y'\vert,\,\,(x,y,y')\in\Bbb R^3\times\Bbb R^3\times\Bbb R^3.
$$
This is the length of the broken path connecting $y$ to $x$ and $x$ to $y'$ which plays the central role
in this paper.

In this paper we denote the {\it convex hull} of the set $F\subset\Bbb R^3$ by $[F]$.

\proclaim{\noindent Theorem 1.1.}
Let $[\overline B\cup\overline B']\cap\partial D=\emptyset$
and $T$ satisfy
$$\displaystyle
T>\min_{x\in\partial D,\,y\in\partial B,\, y'\in\partial B'}\phi(x;y,y').
\tag {1.5}
$$
Then, there exists a $\tau_0>0$ such that, for all $\tau\ge\tau_0$,
$$
\displaystyle
\int_{\Bbb R^3\setminus\overline D}(fv_g-w_fg)dx>0
$$
and the formula
$$\displaystyle
\lim_{\tau\longrightarrow\infty}
\frac{1}{\tau}
\log
\int_{\Bbb R^3\setminus\overline D}(fv_g-w_fg)dx
=-\min_{x\in\partial D,\,y\in\partial B,\, y'\in\partial B'}\phi(x;y,y')
\tag {1.6}
$$
is valid.

\endproclaim
Note that
$$\displaystyle
\min_{x\in\partial D,\, y\in \partial B,\,y'\in \partial B'}\phi(x;y,y')
=\min_{x\in\partial D}\phi(x;p,p')-(\eta+\eta').
\tag {1.7}
$$
See Appendix for the proof of (1.7).
The quantity $\min_{x\in\partial D}\phi(x;p,p')$ coincides with the shortest length
of the broken paths connecting $p$ to a point $q$ on $\partial D$ and $q$ to $p'$,
that is, the {\it first reflection distance} between $p$ and $q$ by $D$.
(1.6) gives an extraction formula of $\min_{x\in\partial D}\phi(x;p,p')$
from $u_f(x,t)$ given at all $x\in B'$ and $t\in]0,\,T[$.
Formula (1.6) gives the method of carrying out calculation processing of the waveform mathematically,
and extracting the first reflection distance.

Define $[p,p']=\{sp+(1-s)p'\,\vert\,0\le s\le 1\}$.  This is the straight line segment connecting
the centers of $B$ and $B'$ and coincides with $[\{p,p'\}]$.  Since both $p$ and $p'$ are in $\Bbb R^3\setminus\overline D$,
$[p,p']\cap\overline D=\emptyset$ if and only if $[p,p']\cap\partial D=\emptyset$.

We know that

$\bullet$  $\min_{x\in\partial D}\phi(x;p,p')\ge\vert p-p'\vert$;

$\bullet$  if $[p,p']\cap\partial D=\emptyset$, then
$\min_{x\in\partial D}\phi(x;p,p')>\vert p-p'\vert$.

Given $c>\vert p-p'\vert$ define
$$\displaystyle
E_c(p,p')=\{x\in\Bbb R^3\,\vert\,\phi(x;p,p')=c\}.
$$
This is a {\it spheroid} with focal points $p$ and $p'$.
It is a compact surface of class $C^{\infty}$.

Since $D$ is contained in the {\it exterior} of spheroid
$E_c(p,p')$ with $c=\min_{x\in\partial D}\phi(x;p,p')$, Theorems
1.1 gives us the largest spheroid with focal points $p$ and $p'$
whose exterior contains $D$ using dynamical bistatic data $u_f$ on
$B'\times\,]0,\,T[$. The appearance of the spheroid in the
enclosure method is new and this is a decisive difference from the
previous enclosure method.

Therefore we obtain the information that there exists a point belonging to $\partial D$ on the spheroid $E_c(p,p')$
with $c=\min_{x\in\partial D}\phi(x;p,p')$ calculated by formula (1.6).
Thus, the next problem is: identify all the points belonging to $\partial D$ on the spheroid.
In order to describe the problem precisely we introduce the following notion.

{\bf\noindent Definition 1.1.} Let $p$ and $p'$ satisfy
$[p,p']\cap\partial D=\emptyset$. Define
$$
\displaystyle
\Lambda_{\partial D}(p,p')
=\{q\in\partial D\,\vert\,
\phi(q;p,p')=\min_{x\in\partial D}\phi(x;p,p')\}.
$$
We call this the {\it first reflector} between $p$ and $p'$.
The points in the first reflector are called the {\it first reflection points} between $p$ and $p'$.
Note that $\Lambda_{\partial D}(p,p')$ can be an infinite set.

One has the expression
$$\displaystyle
\Lambda_{\partial D}(p,p')=\partial D\cap E_c(p,p')
$$
with $c=\min_{x\in\partial D}\phi(x;p,p')$.
Thus the problem becomes: identify all the first reflection points.

Let $\omega\in S^2$.
We denote by $s(\omega;p,p',c)$ the length of the straight line segment connecting $p'$ and
the unique point on $E_c(p,p')\cap\{p'+s\omega\,\vert\, s>0\,\}$.
We have
$$\displaystyle
s(\omega;p,p',c)
=\frac{c^2-\vert p-p'\vert^2}
{2\{c-\omega\cdot(p-p')\}}.
$$
Note that $\omega\cdot(p-p')<c$ since $c>\vert p-p'\vert$.
It is easy to see that the map
$$\displaystyle
S^2\ni\omega\mapsto
p'+s(\omega;p,p',c)\omega
\in\Bbb R^3
$$
is one-to-one and the image coincides with $E_c(p,p')$.

Let $0<\eta'<\inf_{\omega\in S^2}s(\omega;p,p',c)$.
$\overline B'$ is contained in the set
of all $x$ such that $\phi(x;p,p')<c$, that is, the domain enclosed
by $E_c(p,p')$.

The following theorem says that all the first reflection points between $p$ and $p'$
together with the tangent planes
can be extracted from a single set of the bistatic data.
This exceeds the previous enclosure method and suggests that the information which
contained in the bistatic data is quite rich.

\proclaim{\noindent Theorem 1.2.}
Assume that $c=\min_{x\in\partial D}\phi(x;p,p')$ is known.
Let $[\overline B\cup\overline B']\cap\overline D=\emptyset$.
Fix $0<s<\eta'$.
If $T$ satisfies
$$\displaystyle
T>\sup_{\omega\in\,S^2}\min_{x\in\partial D}\phi(x;p,p'+s\omega)-(\eta+\eta'-s),
\tag {1.8}
$$
then, one can extract all $q\in\Lambda_{\partial D}(p,p')$ together with $\nu_q$
from $u_f$ on $B'\times\,]0,\,T[$ with $f=\chi_B$.

\endproclaim

{\bf\noindent Remark 1.1.}  Note that
$$\displaystyle
\sup_{\omega\in\,S^2}\min_{x\in\partial D}\phi(x;p,p'+s\omega)-(\eta+\eta'-s)
=\sup_{\omega\in\,S^2}\min_{x\in\partial D,\,y\in\partial B,\,y'\in\partial B_{\eta'-s}(p'+s\omega)}\phi(x;y,y').
$$
The reason is the same as that of the validity of (1.7).  Thus the constraint on $T$ is reasonable.

Theorem 1.1 is a direct consequence of the following two estimates:
there exist $\mu_j\in\Bbb R$, $C_j>0$ with $j=1,2$ and $\tau_0>0$ which are independent of $\tau$
such that, for all $\tau\ge\tau_0$,
$$\displaystyle
e^{\tau\min_{x\in\partial D,\,y\in\partial B,\,y'\in\partial B'}\phi(x;y,y')}\int_{\Bbb R^3\setminus\overline D}(fv_g-w_fg)dx
\le C_1\tau^{\mu_1}
\tag {1.9}
$$
and
$$\displaystyle
C_2\tau^{\mu_2}
\le e^{\tau\min_{x\in\partial D,\,y\in\partial B,\,y'\in\partial B'}\phi(x;y,y')}\int_{\Bbb R^3\setminus\overline D}(fv_g-w_fg)dx.
\tag {1.10}
$$

The proof of (1.9) proceeds along the same line as the
back-scattering data case ($B=B'$) and is given in Section 2. The point that should be
emphasized in the proof of Theorem 1.1 is (1.10)
which is proved in Subsection 3.1.

When $B'=B$, using the same technique as
done for the {\it sound-hard} obstacle case in \cite{IE3}, we can prove
(1.10) without difficulty.
The technique therein does not depend on the boundary condition, however,
heavily depends on the condition $B=B'$.
In this paper, we take another way. 
It is based on the combination of the {\it maximum
principle} for the modified Helmholtz equation in the domain $\Bbb
R^3\setminus\overline D$ and a {\it reflection} across $\partial
D$.  It heavily depends on the speciality of the homogeneous
Dirichlet boundary condition on $\partial D$. The idea goes back
to the arguments done in the proofs of Theorem 3.6 and Lemma 3.7
in Lax-Phillips \cite{LP}.  Note that, therein, a relationship
between the {\it support function} and the so-called {\it
scattering kernel} for a general sound-soft obstacle has been
established.
They used the arguments to obtain an estimate for the {\it analytic
continuation} of the Fourier transform of the scattering kernel
and then applied the Paley-Weiner theorem.
We refer the reader to \cite{M, PS} for several other results using the scattering kernel.

\subsection{Leading term of the indicator function and its implication}

(1.9) and (1.10) suggest that the following integral as $\tau\longrightarrow\infty$
$$
\displaystyle
e^{\tau\min_{x\in\partial D,\,y\in\partial B,\,y'\in\partial B'}\phi(x;y,y')}\int_{\Bbb R^3\setminus\overline D}(fv_g-w_fg)dx
$$
may behave as some power of $\tau$ multiplied by a positive
constant.  The constant may contain some information about the
geometry of the boundary of the obstacle at the points on
$\partial D$ that attain $\min_{x\in\partial D,\,y\in\partial B,\,y'\in
\partial B'}\phi(x;y,y')$, i.e., the first reflection points between $p$ and $p'$.

If $q\in\Lambda_{\partial D}(p,p')$, then $q\in E_c(p,p')$ with $c=\min_{x\in\partial D}\phi(x;p,p')$
and the two tangent planes at $q$ of $\partial D$ and $E_c(p,p')$ coincide.
We denote by $S_q(\partial D)$ and $S_q(E_c(p,p'))$ the {\it shape operators} (or the {\it Weingarten maps})
at $q$ with respect to $\nu_q$.  Those are symmetric linear operators on the common tangent space at $q$ of $\partial D$ and $E_c(p,p')$.
It is easy to see that $S_q(E_c(p,p'))-S_q(\partial D)\ge 0$ as the quadratic form
on the same tangent space at $q$ (see (4.21)).

Let $Z_f\in H^1(\Bbb R^3\setminus\overline D)$ solve
$$\begin{array}{c}
\displaystyle
(\triangle-\tau^2)Z_f=F_f(x,\tau)\,\,\text{in}\,\Bbb R^3\setminus\overline D,\\
\\
\displaystyle
Z_f=0\,\,\text{on}\,\partial D.
\end{array}
\tag {1.11}
$$
It follows from (1.2) for $h=f$ and (1.4) that $w_f$ has the form
$$\displaystyle
w_f=v_f+\epsilon_f^0+e^{-\tau T}Z_f,
\tag {1.12}
$$
where
$\epsilon_f^0$ satisfies
$$\begin{array}{c}
\displaystyle
(\triangle-\tau^2)\epsilon_f^0=0\,\,\text{in}\,\Bbb R^3\setminus\overline D,\\
\\
\displaystyle
\epsilon_f^0=-v_f\,\,\text{on}\,\partial D.
\end{array}
\tag {1.13}
$$
Note that: since $\text{supp}\,f\cap\overline D=\emptyset$, $v_f$
is smooth in a neighbourhood of $\overline D$ and thus, by
elliptic regularity, we see that $\epsilon_f^0$ is smooth for
$x\in\Bbb R^3\setminus D$.  Moreover, note that
$\epsilon_f^0(x)\longrightarrow 0$ as $\vert
x\vert\longrightarrow\infty$ rapidly and uniformly with respect to $x/\vert
x\vert$. This is a combination of the uniqueness of the weak
solution of (1.13) and a potential theoretic construction of the
solution, see, e.g., \cite{CK, Ms} for the approach and \cite{IK2} for an application
to an inverse problem for the heat equation.

Given $x\in\Bbb R^3$ define $d_{\partial D}(x)=\inf_{y\in\partial D}\vert y-x\vert$.
It is well known that there exists a positive constant
$\delta_0$ such that: given $x\in\overline D$/$x\in \Bbb
R^3\setminus D$ with $d_{\partial D}(x)<2\delta_0$ there
exists a unique $q=q(x)$ be the boundary point on $\partial D$
such that $x=q\mp d_{\partial D}(x)\nu_q$ (\cite{GT}).
One may assume that both $d_{\partial D}(x)$ and $q(x)$ is $C^2$
for $x\in \overline D$ with $d_{\partial D}(x)<2\delta_0$;
$x\in\Bbb R^3\setminus D$ with $d_{\partial D}(x)<2\delta_0$
(Lemma 1 of Appendix in \cite{GT}).
Note that
$\nu_q$ is the unit outer normal to $\partial D$ at $q$.
For $x$ with $d_{\partial D}(x)<2\delta_0$ define
$\displaystyle
x^r=2q(x)-x$.

Before describing our third result, we introduce a restriction on a class of obstacles
which is satisfied with all convex obstacles.

{\bf\noindent Definition 1.2.}
We say that $D$ is {\it admissible}, if
there exist positive constants $C$, $\delta'(\le 2\delta_0)$ and $\tau_0$ such that, for all
$y\in D$ with $d_{\partial D}(y)<\delta'$ and $\tau\ge\tau_0$
$$\displaystyle
\vert\epsilon_f^0(y^r)\vert\le
C\int_B e^{-\tau\vert y-x\vert}dx.
$$

The following theorem gives an answer to the question raised above.

\proclaim{\noindent Theorem 1.3.}
Let $B$ and $B'$ satisfy $[\overline B\cup\overline{B'}]\cap\overline D=\emptyset$.
Let $f=\chi_B$ and $g=\chi_{B'}$.
Let $c=\min_{x\in\partial D}\phi(x;p,p')$.
Let $T$ satisfy (1.5).

Assume that $D$ is admissible and $\partial D$ is $C^3$.
If $\Lambda_{\partial D}(p,p')$ is finite
and for all $q\in\Lambda_{\partial D}(p,p')$
$$\displaystyle
\text{det}\,(S_q(E_c(p,p'))-S_q(\partial D))>0,
\tag {1.14}
$$
then we have
$$\begin{array}{c}
\displaystyle
\lim_{\tau\longrightarrow\infty}
\tau^4e^{\tau\min_{x\in\partial D,\,y\in\partial B,\, y'\in\partial B'}\phi(x;y,y')}
\int_{\Bbb R^3\setminus\overline D}(fv_g-w_fg)dx\\
\\
\displaystyle
=\frac{\pi}{2}\sum_{q\in\Lambda_{\partial D}(p,p')}
\left(\frac{\text{diam}\,B}{2\vert q-p\vert}\right)\cdot
\left(\frac{\text{diam}\,B'}{2\vert q-p'\vert}\right)
\cdot
\frac{1}
{\sqrt{\text{det}\,(S_q(E_c(p,p'))-S_q(\partial D))}}.
\end{array}
\tag {1.15}
$$

\endproclaim

Some remarks are in order.

$\bullet$  The right-hand side of (1.15) is symmetric with respect to the replacement
$p\rightarrow p'$ and $p'\rightarrow p$.
This is a kind of {\it reciprocity}.

$\bullet$  The quantity $\text{det}\,(S_q(E_c(p,p'))-S_q(\partial D))$ expresses
some kind of information about the difference or deviation of the geometry between $\partial D$ and
$E_c(p,p')$ at $q\in\Lambda_{\partial D}(p,p')$.

The following proposition says that Theorem 1.3 can cover convex obstacles.

\proclaim{\noindent Proposition 1.1.}
(i)  If $D$ is convex, then $D$ is admissible and $\Lambda_{\partial D}(p,p')$ consists of a single point.
(ii) Let $q\in\Lambda_{\partial D}(p,p')$ and assume that $\partial D$ is contained in the half space $(x-q)\cdot\nu_q\le 0$,
then (1.14) at $q$ is satisfied.

\endproclaim

For the proof see Appendix.
Thus as a corollary we obtain the following result.

\proclaim{\noindent Corollary 1.1.}
Let $B$ and $B'$ satisfy $[\overline B\cup\overline{B'}]\cap\overline D=\emptyset$.
Let $f=\chi_B$ and $g=\chi_{B'}$.
Let $c=\min_{x\in\partial D}\phi(x;p,p')$.
Let $T$ satisfy (1.5).  If $D$ is convex and $\partial D$ is $C^3$, then (1.15) whose right-hand side
consists of a single term is valid.
\endproclaim

Note that, for the back-scattering case $B=B'$, using the
completely same argument as done in \cite{IE4} in a bounded
domain, we obtain
$$\begin{array}{c}
\displaystyle
\lim_{\tau\longrightarrow\infty}
\tau^4e^{2\tau\text{dist}\,(\partial D,B)}
\int_B(v_f-w_f)dx\\
\\
\displaystyle
=
\frac{\pi}{2}\left(\frac{\text{diam}\,B}{2d_{\partial D}(p)}\right)^2
\sum_{x\in\partial D,\,\vert x-p\vert=d_{\partial D}(p)}
\frac{1}{\sqrt{P_{\partial D}(1/d_{\partial D}(p);x)}},
\end{array}
\tag {1.16}
$$
where $d_{\partial D}(p)=\inf_{x\in\partial D}\vert x-p\vert$;
$\displaystyle
P_{\partial D}(\lambda;q)
=\left(\lambda-k_1(q)\right)
\left(\lambda-k_2(q)\right)$
and $k_1(q)$ and $k_2(q)$ denote the {\it principle curvatures} of $\partial D$ at $q$
with respect to $\nu_q$ (see Appendix in \cite{GT}).  Note that the Gauss curvature $K_{\partial D}(q)$
and mean curvature $H_{\partial D}(q)$ at $q$ with respect to $\nu_q$ are given
by $k_1(q)k_2(q)$ and $(k_1(q)+k_2(q))/2$, respectively.

The assumptions therein are

$\bullet$  $\partial D$ is $C^3$;

$\bullet$  $T>2\text{dist}\,(D,B)$;

$\bullet$  the set of all points $x\in\partial D$ with $\vert x-p\vert=d_{\partial D}(p)$ is {\it finite}
and each point $q$ in the set satisfies
$P_{\partial D}(1/d_{\partial D}(p);q)>0$.

It is not assumed that $D$ is admissible in (1.16) unlike (1.15).

The quantity $P_{\partial D}(1/d_{\partial D}(p);q)$ at $q\in\partial D$ with $\vert q-p\vert=d_{\partial D}(p)$
denotes a `deflection' of the surface $\partial D$ at $q$ from the sphere $\vert x-p\vert=d_{\partial D}(p)$
since we know from, e.g., Proposition 4.2 in this paper that
$\displaystyle
P_{\partial D}(q;1/d_D(p))
=\text{det}\,(S_q(\partial B_{d_{D}(p)}(p))-S_q(\partial D))$,
where $B_{d_D(p)}(p)=\{x\in\Bbb R^3\,\vert\,\vert x-p\vert<d_D(p)\}$.
Thus formula (1.15) of Theorem 1.3 can be considered as an extension of (1.16) to the bistatic data case.
See also Remark 5.2 for a comparison.

After having Theorem 1.3, everyone wishes to extract the geometry
of $\partial D$ at all the first reflection points.  The complete
answer for general obstacle is not known, however, under the
admissibility of $D$, one can obtain the following result.

\proclaim{\noindent Theorem 1.4.}
Let $q\in\Lambda_{\partial D}(p,p')$ be known.
Let $T$ satisfy (1.5).
Let $B$ and $B'$ satisfy $[\overline B\cup\overline B']\cap\overline D=\emptyset$.
If $D$ is admissible and $\partial D$ is $C^3$,
then one can extract
$$
\displaystyle
H_{\partial D}(q)-\frac{S_{q}(\partial D)(\mbox{\boldmath $A$}_q(p)\times\mbox{\boldmath $A$}_q(p'))
\cdot
(\mbox{\boldmath $A$}_q(p)\times\mbox{\boldmath $A$}_q(p'))}{2(1+\mbox{\boldmath $A$}_q(p)\cdot\mbox{\boldmath $A$}_q(p'))}
\tag {1.17}
$$
and
$$\displaystyle
K_{\partial D}(q)
$$
from $u_f$ on $B'\times\,]0,\,T[$ with $f=\chi_B$,
where
$$\displaystyle
\mbox{\boldmath $A$}_q(x)=\frac{q-x}{\vert q-x\vert}.
$$

\endproclaim

Note that $\mbox{\boldmath $A$}_q(p)\times\mbox{\boldmath $A$}_q(p')$ belongs to the tangent space of $\partial D$
at $q$.  For this see Lemma 4.3 in Section 4.

This theorem may suggest the following.

$\bullet$  If one wishes to know the mean curvature at a first reflection point precisely, one
should make the transmitter and the receiver approach as much as possible.
It is because $\mbox{\boldmath $A$}_q(p)\times\mbox{\boldmath $A$}_q(p')$ will disappear approximately at this time
and thus the correction term in (1.17) can be ignored.

$\bullet$  On the other hand, the Gauss curvature at the first reflection point can be extracted regardless of
the position of a transmitter and a receiver at any time
except for the condition $[\overline B\cup\overline B']\cap\overline D=\emptyset$.

As a corollary of Theorems 1.1, 1.2 and 1.4 we have the following result.

\proclaim{\noindent Corollary 1.2.} Assume that $D$ is an open
ball.
Let $T$ satisfy (1.8).
Let $B$ and $B'$ satisfy $[\overline
B\cup\overline B']\cap\overline D=\emptyset$.
Then, one can extract $D$ itself from $u_f$ on
$B'\times\,]0,\,T[$ with $f=\chi_B$.
\endproclaim

The steps to reconstruct an unknown open ball $D$ are as follows.

{\bf Step 1.}  Determine $c=\min_{x\in\partial D}\phi(x;p,p')$ via Theorem 1.1.

{\bf Step 2.}  Determine the unique point $q$ in $\Lambda_{\partial D}(p,p')$ together with $\nu_q$ via Theorem 1.2.

{\bf Step 3.}  Determine $K_{\partial D}(q)$ via Theorem 1.4.

Then the radius and center of $D$ are given by $1/\sqrt{K_{\partial D}(q)}$
and $q-(1/\sqrt{K_{\partial D}(q)})\nu_q$, respectively.

The reconstruction problem of a spherical obstacle also has been considered
in the frequency domain.
For example, see \cite{AMS} which employs a spherical wave as an incident wave
and uses a low frequency limit for the reconstruction.

The four steps described above give a {\it constructive proof}
of a uniqueness theorem in an inverse obstacle problem in the sense that
it does not make use of the {\it uniqueness of the continuation} of the solution
of the governing equation of the wave.
The following uniqueness result employs the bistatic data over a {\it finite time interval}
and itself seems to be new.

\proclaim{\noindent Theorem 1.5.} Let $D_1$ and $D_2$ be open
balls.  Let $u_f^j$ be the solution of (1.1) with $f=\chi_B$ and
$D=D_j$.
Let $T$ satisfy (1.8).
Let $B$ and $B'$ satisfy
$[\overline B\cup\overline B']\cap\overline D_j=\emptyset$ for
$j=1,2$. If $u_f^1=u_f^2$ on $B'\times\,]0,\,T[$, then $D_1=D_2$.

\endproclaim

We refer the readers to \cite {ISA, ISA2, R, R2} for various uniqueness theorems for inverse obstacle problems
for hyperbolic equations over a finite time interval.

Another corollary from Theorem 1.4 is concerned with the
determination of the directions of principle curvatures at a point
on $\partial D$.

Assume that $D$ is convex and $\partial D$ is $C^3$.
From Proposition 1.1 we know that $\Lambda_{\partial D}(p,p')$ consists of a single point.
We denote the point by $q(p,p')$.
We denote by $p(\theta)$ and $p'(\theta)$ the points rotated around the line directed $\nu_q$ at
$q=q(p,p')$ counterclockwise
with rotation angle $\theta\in[0,\,2\pi[$ of $p$ and $p'$.
Thus $p(0)=p$ an $p'(0)=p'$.

Then for all $\theta\in\,[0,\,2\pi[$ we know that $\Lambda_{\partial D}(p(\theta),p'(\theta))$,
$\mbox{\boldmath $A$}_q(p(\theta))\cdot\mbox{\boldmath $A$}_q(p'(\theta))$,
$\vert\mbox{\boldmath $A$}_q(p(\theta))\times\mbox{\boldmath $A$}_q(p'(\theta))\vert$,
$\phi(q(\theta);p(\theta),p'(\theta))$ and $\nu_{q}$ at $q=q(p(\theta),p'(\theta))$
are invariant with respect to $\theta$.

Let $B(\theta)$
denote the
open ball centered at $p(\theta)$ with radius $\eta$
and $B'(\theta)$ the
open ball centered at $p'(\theta)$ with radius $\eta$.
We have $[\overline B(\theta)\cup\overline B'(\theta)]\cap\overline D=\emptyset$
provided $[\overline B\cup\overline B']\cap\overline D=\emptyset$ and $D$ is convex.

Then from (1.17) in Theorem 1.4 applied to $f=f(\theta)$ and $B=B(\theta)$ and $B'=B'(\theta)$
we obtain the function of $\theta$:
$$\displaystyle
\theta\longmapsto
\tilde{H}_{\partial D}(q;p(\theta),p'(\theta))
\equiv
H_{\partial D}(q)-\frac{1}{2}
\sqrt{
\frac{1-\mbox{\boldmath $A$}_q(p)\cdot\mbox{\boldmath $A$}_q(p')}
{1+\mbox{\boldmath $A$}_q(p)\cdot\mbox{\boldmath $A$}_q(p'))}
}
S_{q}(\partial D)(\mbox{\boldmath $V$}(\theta))
\cdot
\mbox{\boldmath $V$}(\theta)
$$
where $\mbox{\boldmath $V$}(\theta)$ denotes the unit vector directed to
$\mbox{\boldmath $A$}_q(p(\theta))\times\mbox{\boldmath $A$}_q(p'(\theta))$.

Now assume that $\mbox{\boldmath $A$}_q(p)\times\mbox{\boldmath $A$}_q(p')\not=0$.
Then, $\mbox{\boldmath $V$}(\theta)$ attains all the tangent vector at
$q$ of $\partial D$ and thus from the behaviour of $\tilde{H}_{\partial D}(q;p(\theta),p'(\theta))$
as a function of $\theta$ one can determine all the directions of principle curvatures say,
$\mbox{\boldmath $V$}(\theta_1)$ and $\mbox{\boldmath $V$}(\theta_2)$ with some $\theta_1$ and $\theta_2$.
Then we have
$$\displaystyle
\frac{\tilde{H}_{\partial D}(q;p(\theta_1),p'(\theta_1))+
\tilde{H}_{\partial D}(q;p(\theta_2),p'(\theta_2))}{2}
=\left\{1-\frac{1}{2}
\sqrt{
\frac{1-\mbox{\boldmath $A$}_q(p)\cdot\mbox{\boldmath $A$}_q(p')}
{1+\mbox{\boldmath $A$}_q(p)\cdot\mbox{\boldmath $A$}_q(p'))}
}\right\}H_{\partial D}(q).
$$
Thus we obtain $H_{\partial D}(q)$.

Summing up, we have obtained the following result.

\proclaim{\noindent Corollary 1.3.}
Let $B$ and $B'$ satisfy $[\overline B\cup\overline B']\cap\overline D=\emptyset$.
Let $T$ satisfy (1.5).
Assume that $D$ is convex and $\partial D$ is $C^3$;
$q=q(p,p')$ is known; $\mbox{\boldmath $A$}_q(p)\times\mbox{\boldmath $A$}_q(p')\not=0$.
Then, one can extract all the directions of principle curvatures, mean and Gauss curvatures, in other words,
the shape operator at $q$ of $\partial D$
from $u_{f(\theta)}$ over $B'(\theta)\times\,]0,\,T[$ for all $\theta\in [0,\,2\pi[$,
where $f(\theta)$ denotes the characteristic function of $B(\theta)$.

\endproclaim

A brief outline of this paper is as follows.
Theorems 1.1 is proved in Sections 2 and 3.
As described above, the key point of the proof is to derive (1.9) and (1.10)
and those are proved in Sections 2 and 3, respectively.

Theorem 1.2 is proved in Subsection 5.1.  The proof contains an
explicit characterization of the first reflector in terms of the
bistatic data.  See Remark 5.1 for the resulted procedure to
determine all the first reflection points.

Theorem 1.3 is proved in Section 4. The key point in the proof of
(1.15) as well as (1.16) is to identify the term which contains
the leading term of the indicator function.  See (4.1) for the
term. We found that the one of two reflection arguments developed
in \cite{LP} works for the purpose. It is based on the reflection
across $\partial D$ and a pointwise estimate of $\epsilon_f^0$
near $\partial D$, that is the use of the admissibility of $D$.
The argument is presented in the proof of Lemma 4.2 in Section 4.
Note that another reflection argument used in the proof of (1.16)
is free from the admissibility assumption, however, can not be
applied to the case when $f\not\equiv g$.

Theorem 1.4 is proved in Subsection 5.2.  The proof is based on an
asymptotic formula which is a consequence of Theorem 1.3 and an
explicit formula of the determinant of the difference of two shape
operators at $q\in\Lambda_{\partial D}(p,p')$ as derived in
Subsection 7.3 of Appendix.

In the final section we give a conclusion of this paper and comments on further problems.

\section{An upper bound of the indicator function}

Define
$$\displaystyle
J(\tau;f,g)
=\int_D(\nabla v_f\cdot\nabla v_g+\tau^2v_fv_g)dx.
\tag {2.1}
$$

We have the following expression of the indicator function.

\proclaim{\noindent Proposition 2.1.}
It follows that
$$\begin{array}{c}
\displaystyle
\int_{\Bbb R^3\setminus\overline D}(fv_g-w_fg)dx
=J(\tau;f,g)
+\int_{\Bbb R^3\setminus\overline D}
(\nabla\epsilon_f^0\cdot\nabla\epsilon_g^0
+\tau^2\epsilon_f^0\epsilon_g^0)dx
\\
\\
\displaystyle
-e^{-\tau T}
\int_{\Bbb R^3\setminus\overline D}(\nabla Z_f\cdot\nabla v_g
+\tau^2 Z_fv_g)dx.
\end{array}
\tag {2.2}
$$

\endproclaim

{\it\noindent Proof.}
From (1.2) and (1.4) we have
$$\displaystyle
\int_{\Bbb R^3\setminus\overline D}(fv_g-w_fg)dx
=\int_{\partial D}\frac{\partial w_f}{\partial\nu}v_gdS
+e^{-\tau T}\int_{\Bbb R^3\setminus\overline D}F_f(x,\tau)v_gdx.
\tag {2.3}
$$
Rewrite
$$\displaystyle
\int_{\partial D}\frac{\partial w_f}{\partial\nu}v_gdS
=\int_{\partial D}\frac{\partial v_f}{\partial\nu}v_gdS
+\int_{\partial D}\frac{\partial\epsilon_f^0}{\partial\nu}v_gdS
+e^{-\tau T}\int_{\partial D}\frac{\partial Z_f}{\partial\nu}v_gdS.
\tag {2.4}
$$
Integration by parts yields
$$\displaystyle
\int_{\partial D}\frac{\partial v_f}{\partial\nu}v_gdS
=J(\tau;f,g).
\tag {2.5}
$$
On the other hand, (1.13) yields
$$\begin{array}{c}
\displaystyle
\int_{\partial D}\frac{\partial\epsilon_f^0}{\partial\nu}v_gdS
=-\int_{\partial D}\frac{\partial\epsilon_f^0}{\partial\nu}\epsilon_g^0dS
=\int_{\Bbb R^3\setminus\overline D}
(\nabla\epsilon_f^0\cdot\nabla\epsilon_g^0
+\tau^2\epsilon_f^0\epsilon_g^0)dx.
\end{array}
\tag {2.6}
$$
Furthermore it follows from (1.2) and (1.11) that
$$\displaystyle
-\int_{\partial D}\frac{\partial Z_f}{\partial\nu}v_gdS
=\int_{\Bbb R^3\setminus\overline D}(\nabla Z_f\cdot\nabla v_g
+\tau^2 Z_fv_g+F_fv_g)dx.
$$
Now from this together with (2.3)-(2.6) we obtain (2.2).

\noindent
$\Box$

\proclaim{\noindent Lemma 2.1.}
As $\tau\longrightarrow\infty$
$$\displaystyle
\left\vert
\int_{\Bbb R^3\setminus\overline D}(\nabla Z_f\cdot\nabla v_g
+\tau^2 Z_fv_g)dx
\right\vert
=O(\tau^{-1}).
\tag {2.7}
$$

\endproclaim

{\it\noindent Proof.}
From (1.2), we obtain
$$\displaystyle
\int_{\Bbb R^3}\left(\vert\nabla v_g\vert^2+
\tau^2\left\vert v_f-\frac{g}{2\tau^2}\right\vert^2\right)dx
=\frac{1}{4\tau^2}\int_B\vert g\vert^2dx.
$$
Since
$$\displaystyle
\left\vert v_g-\frac{g}{2\tau^2}\right\vert^2
\ge
\frac{1}{2}\vert v_g\vert^2
-\frac{\vert g\vert^2}{4\tau^4},
$$
from this we obtain
$$\displaystyle
\frac{1}{2}
\int_{\Bbb R^3}(\vert\nabla v_g\vert^2+\tau^2\vert v_g\vert^2)dx
\le
\frac{1}{2\tau^2}
\int_B\vert g\vert^2dx
$$
and thus
$$\displaystyle
\int_{\Bbb R^3}(\vert\nabla v_g\vert^2+\tau^2\vert v_g\vert^2)dx
\le\frac{1}{\tau^2}\int_B\vert g\vert^2dx.
\tag {2.8}
$$
Similarly it follows from (1.11) that
$$\displaystyle
\int_{\Bbb R^3\setminus\overline D}(\vert\nabla Z_f\vert^2+\tau^2\vert Z_f\vert^2)dx
\le\frac{1}{\tau^2}\int_{\Bbb R^3\setminus\overline D}\vert F_f\vert^2dx.
\tag {2.9}
$$
A combination of (2.8) and (2.9) and the estimate $\Vert F_f\Vert_{L^2(\Bbb R^3\setminus\overline D)}=O(\tau)$ yields (2.7).

\noindent
$\Box$

Thus a combination of (2.2) and (2.7) gives
$$\begin{array}{c}
\displaystyle
\int_{\Bbb R^3\setminus\overline D}(fv_g-w_fg)dx
=J(\tau;f,g)+
\int_{\Bbb R^3\setminus\overline D}
(\nabla\epsilon_f^0\cdot\nabla\epsilon_g^0
+\tau^2\epsilon_f^0\epsilon_g^0)dx
+O(\tau^{-1}e^{-\tau T}).
\end{array}
\tag {2.10}
$$
For the second term in this right-hand side we have the following estimate.

\proclaim{\noindent Lemma 2.2.}
As $\tau\longrightarrow\infty$
$$\displaystyle
\int_{\Bbb R^3\setminus\overline D}
(\vert\nabla \epsilon_f^0\vert^2+\tau^2\vert\epsilon_f^0\vert^2)dx
=O(\tau^2J(\tau;f,f)).
\tag {2.11}
$$

\endproclaim

{\it\noindent Proof.}
It is an application of the trace theorem twice and integration by parts.
More precisely,
choose $\tilde{v}\in H^1(\Bbb R^3\setminus\overline D)$ in such a way that
$$\displaystyle
\Vert \tilde{v}\Vert_{H^1(\Bbb R^3\setminus\overline D)}
\le C\Vert v_f\vert_{\partial D}\Vert_{H^{1/2}(\partial D)},
\tag {2.12}
$$
where $C>0$ and is independent of $v_f$.
Integration by parts (or the weak formulation of (1.13)) yields
$$
\displaystyle
\int_{\partial D}\frac{\partial\epsilon_f^0}{\partial\nu}v_fdS
=-\int_{\Bbb R^3\setminus\overline D}
(\nabla\epsilon_f^0\cdot\nabla\tilde{v}+\tau^2\epsilon_f^0\tilde{v})dx.
\tag {2.13}
$$
A combination (2.12) and (2.13) gives
$$
\displaystyle
\left\vert\int_{\partial D}\frac{\partial\epsilon_f^0}{\partial\nu}v_fdS\right\vert
\le C(\Vert\nabla\epsilon_f^0\Vert_{L^2(\Bbb R^3\setminus\overline D)}
+\tau^2\Vert\epsilon_f^0\Vert_{L^2(\Bbb R^3\setminus\overline D)})
\Vert v_f\vert_{\partial D}\Vert_{H^{1/2}(\partial D)}.
\tag {2.14}
$$
Since
$$\displaystyle
\max{(\Vert\nabla\epsilon_f^0\Vert_{L^2(\Bbb R^3\setminus\overline D)},\tau\Vert\epsilon_f^0\Vert_{L^2(\Bbb R^3\setminus\overline D)})}
\le
\left(\int_{\Bbb R^3\setminus\overline D}
(\vert\nabla \epsilon_f^0\vert^2+\tau^2\vert\epsilon_f^0\vert^2)dx\right)^{1/2},
$$
it follows from (2.14) that
$$\displaystyle
\left\vert\int_{\partial D}\frac{\partial\epsilon_f^0}{\partial\nu}v_fdS\right\vert
\le C(1+\tau)\left(\int_{\Bbb R^3\setminus\overline D}
(\vert\nabla \epsilon_f^0\vert^2+\tau^2\vert\epsilon_f^0\vert^2)dx\right)^{1/2}\Vert v_f\vert_{\partial D}\Vert_{H^{1/2}(\partial D)}.
$$
From this together with (2.6) for $f=g$, we obtain
$$\displaystyle
\int_{\Bbb R^3\setminus\overline D}
(\vert\nabla \epsilon_f^0\vert^2+\tau^2\vert\epsilon_f^0\vert^2)dx
\le C^2(1+\tau)^2
\Vert v_f\vert_{\partial D}\Vert_{H^{1/2}(\partial D)}^2.
\tag {2.15}
$$
By the trace theorem, we have
$$\displaystyle
\Vert v_f\vert_{\partial D}\Vert_{H^{1/2}(\partial D)}^2
\le C'
(\Vert\nabla v_f\Vert_{L^2(D)}^2+\Vert v_f\Vert_{L^2(D)}^2),
$$
where $C'>0$ is independent of $v_f$.
Now from this together with (2.15) and the trivial estimates
$$\displaystyle
\max{(\Vert\nabla v_f\Vert_{L^2(D)}^2, \tau^2\Vert v_f\Vert_{L^2(D)}^2})\le J(\tau;f,f),
$$
yields (2.11).

\noindent
$\Box$

Therefore (1.9) with $\mu_1=2$ follows from (2.10) and (2.11) together with the following estimate.

\proclaim{\noindent Lemma 2.3.}
We have, as$\tau\longrightarrow\infty$,
$$\displaystyle
J(\tau;f,g)=O(\tau e^{-\tau\min_{x\in\partial D,\,y\in\partial B,\,y'\in\partial B'}\phi(x;y,y')}).
\tag {2.16}
$$

\endproclaim

{\it\noindent Proof.}
It follows from (1.3) and (2.5) that
$$\displaystyle
J(\tau;f,g)
=\left(\frac{1}{4\pi}\right)^2
\int_{\partial D}dS_x
\int_{B\times B'}k_{\tau}(x,y,y')dydy',
\tag {2.17}
$$
where
$$\displaystyle
k_{\tau}(x,y,y')
=\left(\frac{1}{\vert x-y\vert}+\tau\right)
\frac{(y-x)\cdot\nu_x}{\vert x-y\vert^2\vert x-y'\vert}e^{-\tau(\vert x-y\vert+\vert x-y'\vert)},
\,\,(x,y,y')\in\partial D\times B\times B'.
\tag {2.18}
$$
Since we have
$$\displaystyle
\inf_{x\in\partial D,\,y\in\,B,\,y'\in\,B'}\phi(x;y,y')
=\min_{x\in\partial D,\,y\in\partial B,\,y'\in\partial B'}\phi(x;y,y')
$$
and $\overline B\cap\overline D=\overline B'\cap\overline D=\emptyset$, from (2.17) and (2.18)
we obtain (2.16).

\noindent
$\Box$

\section{A lower bound of the indicator function}

\subsection{A reduction to a convex obstacle and the proof of (1.10).}

Rewriting the second term in the right-hand side of (2.10) with (2.6),
one has
$$\begin{array}{c}
\displaystyle
\displaystyle
\int_{\Bbb R^3\setminus\overline D}(fv_g-w_fg)dx
=J(\tau;f,g)
+\int_{\partial D}\frac{\partial\epsilon_f^0}{\partial\nu}v_gdS
+O(\tau^{-1}e^{-\tau T}).
\end{array}
\tag {3.1}
$$
In general, we do not know the signature of the second term of the right-hand side of (3.1),
however, we know that the function under integral is nonnegative at a special
point on $\partial D$ by virtue of
the following lemma which is an application of the maximum principle for differential operator $\triangle-\tau^2$
and a reflection argument in \cite{LP}.  It corresponds to Lemma 3.7 in \cite{LP} in which $v_f$ in (1.13) is
replaced with $-e^{-\tau x\cdot\omega}$ for a $\omega\in S^2$.

\proclaim{\noindent Lemma 3.1.}
Let $q\in\partial D$ be a point of support of $D$, i.e., $D$ is contained in the half-space $x\cdot n_q<q\cdot n_q$.
We have
$$\displaystyle
\frac{\partial\epsilon_f^0}{\partial\nu}(q)\ge 0.
\tag {3.2}
$$
\endproclaim

{\it\noindent Proof.}
First,  we prove that, for all $x\in\Bbb R^3\setminus\overline D$,
$$\displaystyle
\epsilon_f^0(x)\ge -v_f(x).
\tag {3.3}
$$
Define $w_f^0=\epsilon_f^0+v_f\in H^1(\Bbb R^3\setminus\overline D)$.
It holds that
$$\begin{array}{c}
\displaystyle
(\triangle-\tau^2)w_f^0=-f\,\,\text{in}\,\Bbb R^3\setminus\overline D,\\
\\
\displaystyle
w_f^0=0\,\,\text{on}\,\partial D.
\end{array}
$$
Since $(\triangle-\tau^2)w_f^0\le 0$ in $\Bbb
R^3\setminus\overline D$, from the weak maximum principle for
operators of divergence form (Theorem 8.1 in \cite{GT}), we obtain
$$\displaystyle
\inf_{(\Bbb R^3\setminus\overline D)\cap B_R}w_f^0
\ge\min\,(\inf_{S_R}w_f^0,0).
\tag {3.4}
$$
where $B_R$ denotes an arbitrary open ball centered at the origin
such that $\overline D\cup\text{supp}\,f \subset B_R$ and
$S_R=\partial B_R$. Since $\epsilon_f^0$ decays as $\vert
x\vert\longrightarrow\infty$ uniformly with respect to $x/\vert
x\vert$, we have $\inf_{S_R}w_f^0\longrightarrow 0$ as
$R\longrightarrow\infty$, and thus, from (3.4) we obtain
$w_f^0(x)\ge 0$ for all $x\in\Bbb R^3\setminus\overline D$. This
completes the proof of (3.3).

The equality in (3.3) holds for $x=q$.  This implies the following inequality for the normal derivatives:
$$\displaystyle
\frac{\partial\epsilon_f^0}{\partial\nu}(q)\ge
-\frac{\partial v_f}{\partial\nu}(q).
\tag {3.5}
$$

Second, we prove that, for all points that satisfy $x\cdot n_q\ge q\cdot n_q$,
$$\displaystyle
\epsilon_f^0(x)\ge -v_f(x'),
\tag {3.6}
$$
where $x'$ is the image of $x$ under reflection across the plane $x\cdot n_q=q\cdot n_q$.
Since $x=x'$ on the plane $x\cdot n_q=q\cdot n_q$, (3.3) shows that (3.6) is satisfied there.
Note that also $v_f'(x)\equiv v_f(x')$ satisfies $(\triangle-\tau^2)v_f'=-f(x')\le 0$.
Applying the weak maximum principle to
$(w_f^0)'\equiv\epsilon_f^0+v_f'$
in the half-space $x\cdot n_q>q\cdot n_q$, one obtains as before (3.6) holds throughout the half-space.
This completes the proof of (3.6).

Since the equality in (3.6) holds for $x=q$, it follows as before that
$$\displaystyle
\frac{\partial\epsilon_f^0}{\partial\nu}(q)
\ge -\frac{\partial}{\partial\nu}\left\{v_f(x')\right\}\vert_{x=q}
=\frac{\partial v_f}{\partial\nu}(q).
\tag {3.7}
$$
Now a combination of (3.5) and (3.7) yields (3.2).

\noindent
$\Box$

The following lemma is an easy consequence
of the $C^2$-regularity of $\partial D$ and thus the proof is omitted.

\proclaim{\noindent Lemma 3.2.}
Let $q\in\Lambda_{\partial D}(p,p')$.
Then, there exists an open ball $\tilde{D}$ contained in $D$ such that
$q\in\Lambda_{\partial\tilde{D}}(p,p')$ and thus $\min_{x\in\partial\tilde{D}}\phi(x;p,p')=\min_{x\in\partial D}\phi(x;p,p')$.
\endproclaim

Let $\tilde{u}=\tilde{u}_f$ denote the weak solution of the following initial boundary value problem:
$$\begin{array}{c}
\displaystyle
\partial_t^2\tilde{u}-\triangle\tilde{u}=0\,\,\text{in}\,(\Bbb R^3\setminus\overline{\tilde{D}})\times\,]0,\,T[,\\
\\
\displaystyle
\tilde{u}(x,0)=0\,\,\text{in}\,\Bbb R^3\setminus\overline{\tilde{D}},\\
\\
\displaystyle
\partial_t\tilde{u}(x,0)=f(x)\,\,\text{in}\,\Bbb R^3\setminus\overline{\tilde{D}},\\
\\
\displaystyle
\tilde{u}=0\,\,\text{on}\,\partial\tilde{D}\times\,]0,\,T[.
\end{array}
\tag {3.8}
$$
Define
$$\displaystyle
\tilde{w}_f(x,\tau)
=\int_0^Te^{-\tau T}\tilde{u}(x,t)dt,\,\,x\in\Bbb R^3\setminus\overline{\tilde{D}},\,\tau>0.
$$

\proclaim{\noindent Lemma 3.3.}
We have
$$\displaystyle
\int_{\Bbb R^3\setminus\overline D}(fv_g-w_fg)dx
\ge\int_{\Bbb R^3\setminus\overline{\tilde{D}}}(fv_g-\tilde{w}_fg)dx+O(\tau^{-1}e^{-\tau T}).
\tag {3.9}
$$
\endproclaim

{\it\noindent Proof.}
Let $\tilde{Z}_f\in H^1(\Bbb R^3\setminus\overline{\tilde{D}})$ solve
$$\begin{array}{c}
\displaystyle
(\triangle-\tau^2)\tilde{Z}_f=\tilde{F}_f(x,\tau)\,\,\text{in}\,\Bbb R^3\setminus\overline{\tilde{D}},\\
\\
\displaystyle
\tilde{Z}_f=0\,\,\text{on}\,\partial\tilde{D},
\end{array}
$$
where
$$\displaystyle
\tilde{F}_f(x,\tau)
=\partial_t\tilde{u}_f(x,T)+\tau\tilde{u}_f(x,T),\,\,x\in\Bbb R^3\setminus\overline{\tilde{D}}.
$$
Similar to $Z_f$ which is the solution of (2.2), we have $\Vert\tilde{Z}_f\Vert_{L^2(\Bbb R^3\setminus\overline{\tilde{D}})}=O(\tau^{-1})$.
And similar to (1.12) for $w_f$, $\tilde{w}_f$ has the form
$$\displaystyle
\tilde{w}_f=v_f+\tilde{\epsilon}_f^0+e^{-\tau T}\tilde{Z}_f,
$$
where
$\tilde{\epsilon}_f^0$ satisfies
$$\begin{array}{c}
\displaystyle
(\triangle-\tau^2)\tilde{\epsilon}_f^0=0\,\,\text{in}\,\Bbb R^3\setminus\overline{\tilde{D}},\\
\\
\displaystyle
\tilde{\epsilon}_f^0=-v_f\,\,\text{on}\,\partial\tilde{D}.
\end{array}
$$
Thus, we have
$$\displaystyle
\tilde{w}_f-w_f=(\tilde{\epsilon}_f^0-\epsilon_f^0)+e^{-\tau T}(\tilde{Z}_f-Z_f)\,\,\text{in}\,\Bbb R^3\setminus\overline D.
$$
Since $\tilde{v_f}\ge 0$ on $\partial\tilde{D}$ and
$\tilde{\epsilon}_f^0(x)\longrightarrow 0$ as
$\vert x\vert\longrightarrow\infty$, by the maximum principle for the
modified Helmholtz equation in $\Bbb
R^3\setminus\overline{\tilde{D}}$, we have
$-\tilde{\epsilon}_f^0\le v_f$ in $\Bbb R^3\setminus\tilde{D}$ and
thus $-\tilde{\epsilon}_f^0\le-\epsilon_f^0$ on $\partial D$.
Again by the maximum principle for the modified Helmholtz equation
in $\Bbb R^3\setminus\overline D$, we obtain
$-\tilde{\epsilon}_f^0\le -\epsilon_f^0$ in $\Bbb R^3\setminus D$.
Therefore we obtain
$$\displaystyle
\tilde{w}_f-w_f\ge e^{-\tau T}(\tilde{Z}_f-Z_f)\,\,\text{in}\,\Bbb R^3\setminus\overline D.
\tag {3.10}
$$
Since both $\text{supp}\,g$ and $\text{supp}\,f$ are contained in $\Bbb R^3\setminus\overline D$ and thus
combining this with (3.10), we obtain
$$\displaystyle
\int_{\Bbb R^3\setminus\overline D}(fv_g-w_fg)dx
-\int_{\Bbb R^3\setminus\overline{\tilde{D}}}(fv_g-\tilde{w}_fg)dx
\ge e^{-\tau T}\int_{\Bbb R^3\setminus\overline D}(\tilde{Z}_f-Z_f)gdx.
$$
From the $L^2$-bounds for $Z_f$ and $\tilde{Z}_f$, we see that this right-hand side has the bound $O(\tau^{-1}e^{-\tau T})$.

\noindent
$\Box$

Since $\tilde{D}$ is convex, every point $q\in\partial\tilde{D}$ is a point of support of $\tilde{D}$ and thus,
from (3.1), (3.2) and (3.9) we obtain
$$\begin{array}{c}
\displaystyle
\displaystyle
\int_{\Bbb R^3\setminus\overline D}(fv_g-w_fg)dx
\ge \tilde{J}(\tau;f,g)
+O(\tau^{-1}e^{-\tau T}),
\end{array}
\tag {3.11}
$$
where
$$\displaystyle
\tilde{J}(\tau;f,g)
=\int_{\tilde{D}}(\nabla v_f\cdot\nabla v_g+\tau^2v_fv_g)dx.
$$
Now everything is reduced to give a lower estimate for $\tilde{J}(\tau;f,g)$ as $\tau\longrightarrow\infty$.
For this and the future use of it in the sound-hard obstacle case we give the estimate for $J(\tau;f,g)$ for
general $D$.

In the following lemma we do not assume that $D$ is convex.

\proclaim{\noindent Lemma 3.4.}
There exist positive constants $C$, $\mu$ and $\tau_0$ such that,
for all $\tau\ge\tau_0$,
$$\displaystyle
\tau^{2+\mu}e^{\tau\min_{x\in\partial D}\phi(x;p,p')}
e^{-\tau(\eta+\eta')}J(\tau;f,g)
\ge C.
\tag {3.12}
$$

\endproclaim

We give the proof of this lemma in the next subsection.

It follows from (3.11) and (3.12) for $D=\tilde{D}$ that there exist positive constants $C'$ and $\tau_0'>0$ such that,
for all $\tau\ge\tau_0'$,
$$
\displaystyle
\displaystyle
\tau^{2+\mu}e^{\tau\min_{x\in\partial D,\,y\in\partial B,y'\in\partial B'}\phi(x;y,y')}
\int_{\Bbb R^3\setminus\overline D}(fv_g-w_fg)dx\ge C'
$$
provided $T$ satisfies (1.5).  This completes the proof of (1.10).

\subsection{Proof of Lemma 3.4.}

In this subsection we never assume that $D$ is convex.
Let $A(x,\tau)$ be an arbitrary positive function of $x\in\partial D$ with parameter $\tau>0$.
For another function $B(x,\tau)$,
in the following,  $B(x,\tau)=O(A(x,\tau))$ as $\tau\longrightarrow\infty$
and uniformly with respect to $x\in\partial D$ means
that there exist positive constants $\tau_0$ and $C$ independent of $x\in\partial D$ such that,
for all $\tau\ge\tau_0$ and $x\in\partial D$ we have
$\displaystyle
\vert B(x,\tau)\vert\le CA(x,\tau)$.

The proof of Lemma 3.4 starts with having the following expression.

\proclaim{\noindent Lemma 3.5.}
There exists a positive constant $C$ such that, as $\tau\longrightarrow\infty$
$$\begin{array}{c}
\displaystyle
J(\tau;f,g)\\
\\
\displaystyle
=\frac{1}{4\tau^3}\left(\eta-\frac{1}{\tau}\right)\left(\eta'-\frac{1}{\tau}\right)
\int_{\partial D}
\frac{(p-x)\cdot\nu_x}{\vert x-p\vert^2\vert x-p'\vert}
\left(1+\frac{1}{\tau\vert x-p\vert}\right)
e^{-\tau(d_B(x)+d_{B'}(x))}dS_x\\
\\
\displaystyle
+O(\tau^{-1}e^{-\tau\inf_{x\in\partial D}(d_B(x)+d_{B'}(x))}
(e^{-C\tau\min_{x\in\partial D}d_B(x)}+e^{-C\tau\min_{x\in\partial D}d_{B'}(x)})).
\end{array}
\tag {3.13}
$$

\endproclaim

{\it\noindent Proof.}
By \cite{IE4}, we have, as
$\tau\longrightarrow\infty$ and uniformly with respect to
$x\in\partial D$,
$$\displaystyle
v_g(x)
=\frac{1}{2\tau}
\left((J_1^+)_g(x,\tau)+O(e^{-\tau d_{B'}(x)(1+C)})\right)
$$
and
$$\displaystyle
\nabla v_f(x)
=\frac{1}{2}
\left((J_2^+)_f(x,\tau)\frac{p-x}{\vert x-p\vert}
+O(e^{-\tau d_B(x)(1+C)})\right),
$$
where $C$ is a positive constant,
$$\begin{array}{c}
\displaystyle
(J_2^+)_f(x,\tau)
=\frac{e^{-\tau d_B(x)}}{\tau\vert x-p\vert}
\left(\eta-\frac{1}{\tau}\right)
\left(1+\frac{1}{\tau\vert x-p\vert}\right)\\\\
\displaystyle
+\frac{e^{-\tau d_B(x)\sqrt{1+(2\eta/d_B(x))}}}{\tau^2\vert x-p\vert^2}
\left(d_B(x)\sqrt{1+\frac{2\eta}{d_B(x)}}
+\frac{1}{\tau}\right)
\end{array}
$$
and
$$\begin{array}{c}
\displaystyle
(J_1^+)_g(x,\tau)
=\frac{e^{-\tau d_{B'}(x)}}
{\tau\vert x-p'\vert}
\left(\eta'-\frac{1}{\tau}\right)
+\frac{e^{-\tau d_{B'}(x)\sqrt{1+(2\eta'/d_{B'}(x))}}}{\tau^2\vert x-p'\vert}.
\end{array}
$$
From these and (2.5) we obtain
$$\begin{array}{c}
\displaystyle
J(\tau;f,g)
=\frac{1}{4\tau}
\int_{\partial D}(J_2^+)_f(x,\tau)\frac{(p-x)\cdot\nu_x}{\vert x-p\vert}
(J_1^+)_g(x,\tau)dS_x\\
\\
\displaystyle
+O(\tau^{-1}e^{-\tau\inf_{x\in\partial D}(d_B(x)+d_{B'}(x))}
(e^{-C\tau\min_{x\in\partial D}d_B(x)}+e^{-C\tau\min_{x\in\partial D}d_{B'}(x)})).
\end{array}
$$
This yields (3.13) since we have as $\tau\longrightarrow\infty$ and uniformly with respect to $x\in\partial D$,
$$\begin{array}{c}
\displaystyle
(J_1^+)_g(x,\tau)(J_2^+)_f(x,\tau)
=\frac{e^{-\tau(d_B(x)+d_{B'}(x))}}{\tau^2\vert x-p\vert\vert x-p'\vert}
\left(\eta-\frac{1}{\tau}\right)\left(\eta'-\frac{1}{\tau}\right)
\left(1+\frac{1}{\tau\vert x-p\vert}\right)\\
\\
\displaystyle
+O(\tau^{-3}e^{-\tau(d_B(x)+d_{B'}(x))}
(e^{-\tau d_{B'}(x)C}
+e^{-\tau d_B(x)C}))
+O(\tau^{-4}e^{-\tau(d_B(x)+d_{B'}(x))(1+C)}).
\end{array}
$$

\noindent
$\Box$

Now we give a lower estimate of $J(\tau;f,g)$ as $\tau\longrightarrow\infty$ by using (3.13).

Define
$$\displaystyle
I_m(\tau)
=\int_{\partial D}
\frac{(p-x)\cdot\nu_x}{\vert x-p\vert^m\vert x-p'\vert}
e^{-\tau\phi(x;p,p')}dS_x,
$$
where $m=2,3$.

Since $d_B(x)=\vert x-p\vert-\eta$ for $x\in\Bbb R^3\setminus B$ and $d_{B'}(x)=\vert x-p'\vert-\eta'$
for $x\in\Bbb R^3\setminus B'$, it follows from (3.13) that
$$\begin{array}{c}
\displaystyle
e^{-\tau(\eta+\eta')}J(\tau;f,g)
=\frac{1}{4\tau^3}\left(\eta-\frac{1}{\tau}\right)\left(\eta'-\frac{1}{\tau}\right)
\left(I_2(\tau)+\frac{1}{\tau}I_3(\tau)\right)\\
\\
\displaystyle
+O(\tau^{-1}e^{-\tau\min_{x\in\partial D}\phi(x;p,p')}
(e^{-C\tau\min_{x\in\partial D}d_B(x)}+e^{-C\tau\min_{x\in\partial D}d_{B'}(x))}).
\end{array}
\tag {3.14}
$$
Since
$$\displaystyle
\nabla\cdot\left\{
\frac{(p-x)}{\vert x-p\vert^m\vert x-p'\vert}\right\}
=\frac{(m-3)}{\vert x-p\vert^m\vert x-p'\vert}
+\frac{(p-x)\cdot(p'-x)}{\vert x-p\vert^m\vert x-p'\vert^3}
$$
and
$$\displaystyle
(p-x)\cdot\nabla\phi(x;p,p')
=-\left(\vert x-p\vert+\frac{(p-x)\cdot(p'-x)}{\vert x-p'\vert}\right),
$$
we have
$$\displaystyle
\begin{array}{c}
\displaystyle
I_m(\tau)
=(m-3)\int_D\frac{e^{-\tau\phi(x;p,p')}}{\vert x-p\vert^m\vert x-p'\vert}dx
+\int_D\frac{(p-x)\cdot(p'-x)}{\vert x-p\vert^m\vert x-p'\vert^3}
e^{-\tau\phi(x;p,p')}dx\\
\\
\displaystyle
-\tau\int_D
\frac{(p-x)\cdot\nabla\phi(x;p,p')}{\vert x-p\vert^m\vert x-p'\vert}
e^{-\tau\phi(x;p,p')}dx\\
\\
\displaystyle
=(m-3)\int_D\frac{e^{-\tau\phi(x;p,p')}}{\vert x-p\vert^m\vert x-p'\vert}dx
+\int_D\frac{(p-x)\cdot(p'-x)}{\vert x-p\vert^m\vert x-p'\vert^3}
e^{-\tau\phi(x;p,p')}dx\\
\\
\displaystyle
+\tau
\int_D
\frac{e^{-\tau\phi(x;p,p')}}{\vert x-p\vert^{m-1}\vert x-p'\vert}dx
+\tau
\int_D\frac{(p-x)\cdot(p'-x)}{\vert x-p\vert^m\vert x-p'\vert^2}
e^{-\tau\phi(x;p,p')}dx.
\end{array}
$$
This yields
$$\begin{array}{c}
\displaystyle
I_2(\tau)+\frac{1}{\tau}I_3(\tau)
=
\tau\int_D
\left\{1+\frac{(p-x)\cdot(p'-x)}{\vert x-p\vert\vert x-p'\vert}\right\}
\frac{e^{-\tau\phi(x;p,p')}}{\vert x-p\vert\vert x-p'\vert}dx\\
\\
\displaystyle
+\int_D\left(\frac{1}{\vert x-p\vert}+\frac{1}{\vert x-p'\vert}\right)
\frac{(p-x)\cdot(p'-x)}{\vert x-p\vert^2\vert x-p'\vert^2}e^{-\tau\phi(x;p,p')}dx\\
\\
\displaystyle
+
\frac{1}{\tau}\int_D\frac{(p-x)\cdot(p'-x)}{\vert x-p\vert^3\vert x-p'\vert^3}e^{-\tau\phi(x;p,p')}dx.
\end{array}
\tag {3.15}
$$

\proclaim{\noindent Lemma 3.6.} Let $D$ be an arbitrary nonempty
bounded open set. If $p$ and $p'$ be arbitrary points in $\Bbb
R^3\setminus\overline D$ such that $[p,p']\cap\overline D=\emptyset$, then
$$\displaystyle
C_D(p,p')\equiv\inf_{x\in D}
\left\{1+\frac{(p-x)\cdot(p'-x)}{\vert p-x\vert\vert p'-x\vert}\right\}>0.
\tag {3.16}
$$

\endproclaim

{\it\noindent Proof.}
It is easy to see that $0\le C_D(p,p')\le 2$.
Assume that $C_D(p,p')=0$.  Since we have the identity
$$\displaystyle
1+\frac{(p-x)\cdot(p'-x)}{\vert p-x\vert\vert p'-x\vert}
=\frac{1}{2}
\left\vert
\frac{p-x}{\vert p-x\vert}+\frac{p'-x}{\vert p'-x\vert}\right\vert^2,\,x\not=p,p',
$$
there exists a sequence $\{x_n\}$ in $D$ such that
$$
\displaystyle
\left\vert
\frac{p-x_n}{\vert p-x_n\vert}+\frac{p'-x_n}{\vert p'-x_n\vert}\right\vert^2\longrightarrow 0.
\tag {3.17}
$$
Since $\overline D$ is compact, choosing a subsequence of $\{x_n\}$ if necessary,
one may assume that $\{x_n\}$ converges to a point $y\in\overline D$.  Since $p,p'\not=y$ by assumption,
it follows from (3.17) that
$$\displaystyle
\frac{p-y}{\vert p-y\vert}+\frac{p'-y}{\vert p'-y\vert}=0.
$$
This gives $y\in\,[p,p']$ and thus $[p,p']\cap\overline D\not=\emptyset$.
This is a contradiction.

\noindent
$\Box$

A combination of (3.15) and (3.16) gives
$$\displaystyle
I_2(\tau)+\frac{1}{\tau}I_3(\tau)
\ge \tau K(\tau)\int_D e^{-\tau\phi(x;p,p')}dx,
\tag {3.18}
$$
where
$$\displaystyle
K(\tau)
=\frac{C_D(p,p')}{d_{\partial D}(p)d_{\partial D}(p')}
-\left(\frac{1}{d_{\partial D}(p)}+\frac{1}{d_{\partial D}(p')}\right)
\frac{1}{d_{\partial D}(p)d_{\partial D}(p')}\frac{1}{\tau}
-\frac{1}{d_{\partial D}(p)^2d_{\partial D}(p')^2}\frac{1}{\tau^2}.
$$

\proclaim{\noindent Lemma 3.7.}
Let $p, p'\in\Bbb R^3$.
Then, there exists a number $\mu$ such that
$$\displaystyle
\liminf_{\tau\longrightarrow\infty}\tau^{\mu}e^{\tau\min_{x\in\partial D}\phi(x;p,p')}
\int_D e^{-\tau\phi(x;p,p')}dx>0.
\tag {3.19}
$$

\endproclaim

{\it\noindent Proof.}
Let $x_0\in\partial D$ be a point such that $\phi(x_0;p,p')=\min_{x\in\partial D}\phi(x;p,p')$.
Since $\vert p-x\vert\le\vert p-x_0\vert+\vert x_0-x\vert$ and $\vert p'-x\vert\le\vert p'-x_0\vert+\vert x_0-x\vert$,
we have
$$\displaystyle
\phi(x;p,p')
\le\phi(x_0;p,p')+2\vert x_0-x\vert,\,\,\forall x\in\Bbb R^3.
$$
This gives
$$\displaystyle
e^{\tau\min_{x\in\partial D}\phi(x;p,p')}\int_D e^{-\tau\phi(x;p,p')}dx
\ge
\int_D e^{-2\tau\vert x-x_0\vert}dx.
$$
In \cite{IK2} we have already known that
$$\displaystyle
\liminf_{\tau\longrightarrow\infty}\tau^{3}
\int_D e^{-2\tau\vert x-x_0\vert}dx>0.
$$
Thus (3.19) is valid for $\mu=3$.

\noindent
$\Box$

Now it follows from (3.14), (3.18) and (3.19) that (3.12) is valid.

\section{Asymptotic behaviour of the indicator function}

First we claim that
$$\begin{array}{c}
\displaystyle
\int_{\Bbb R^3\setminus\overline D}(fv_g-w_fg)dx
=2J(\tau;f,g)(1+O(\tau^{-1/2}))
+O(\tau^{-1}e^{-\tau T}).
\end{array}
\tag {4.1}
$$
This is a consequence of the following asymptotic formula and (2.10).

\proclaim{\noindent Proposition 4.1.}
As $\tau\longrightarrow\infty$,
$$\displaystyle
\int_{\Bbb R^3\setminus\overline D}
(\nabla\epsilon_f^0\cdot\nabla\epsilon_g^0+\tau^2
\epsilon_f^0\epsilon_g^0)dx
=J(\tau;f,g)(1+O(\tau^{-1/2})).
$$
\endproclaim

Proposition 4.1 is a direct consequence of the following lemmas.

\proclaim{\noindent Lemma 4.1.} Let $D$ be an arbitrary nonempty
bounded open set. If $B$ and $B'$ be arbitrary open balls such
that $[\overline B\cup\overline B']\cap\overline D=\emptyset$, then
there exist positive constants $C$ and $\tau_0$ independent of $\tau$ such that, for all $\tau\ge\tau_0$,
$J(\tau;f,g)>0$ and
$$\displaystyle
J(\tau;f,g)\ge C\tau^2\int_Ddx\int_{B\times B'} e^{-\tau\phi(x;y,y')}dydy'.
$$

\endproclaim

{\it\noindent Proof.}
From (1.3) and (2.1) we have
$$\displaystyle
J(\tau;f,g)
=\left(\frac{1}{4\pi}\right)^2
\int_Ddx\int_{B\times B'}dydy'\frac{K_{\tau}(x,y,y')}{\vert x-y\vert\vert x-y'\vert}e^{-\tau\phi(x;y,y')},
\tag {4.2}
$$
where
$$\displaystyle
K_{\tau}(x,y,y')=\left(\frac{1}{\vert x-y\vert}+\tau\right)
\left(\frac{1}{\vert x-y'\vert}+\tau\right)
\frac{(y-x)}{\vert x-y\vert}\cdot\frac{(y'-x)}{\vert x-y'\vert}
+\tau^2.
$$
Since $\overline B\cap\overline D=\overline B'\cap\overline D=\emptyset$, there exist positive constants
$C_1$ and $C_2$ such that
$$\displaystyle
K_{\tau}(x,y,y')
\ge
\tau^2\left(1+\frac{(y-x)}{\vert x-y\vert}\cdot\frac{(y'-x)}{\vert x-y'\vert}\right)
-C_1\tau-C_2,\,\,\forall (x,y,y')\in D\times B\times B'.
\tag {4.3}
$$
Here we claim the following estimate:
$$\displaystyle
C_D(B,B')\equiv\inf_{(x,y,y')\in D\times B\times B'}
\left(1+\frac{y-x}{\vert y-x\vert}\cdot\frac{y'-x}{\vert y'-x\vert}\right)>0.
\tag {4.4}
$$
This is proved as follows.

It is easy to see that $0\le C_D(B,B')\le 2$.
Assume that $C_D(B,B')=0$.  Since we have the identity
$$\displaystyle
1+\frac{y-x}{\vert y-x\vert}\cdot\frac{y'-x}{\vert y'-x\vert}
=\frac{1}{2}
\left\vert
\frac{y-x}{\vert y-x\vert}+\frac{y'-x}{\vert y'-x\vert}\right\vert^2,\,x\not=y,y',
$$
there exist sequences $\{x_n\}$ in $D$, $\{y_n\}$ in $B$ and $\{y_n'\}$ in $B'$ such that
$$
\displaystyle
\left\vert
\frac{y_n-x_n}{\vert y_n-x_n\vert}+\frac{y_n'-x_n}{\vert y_n'-x_n\vert}\right\vert^2\longrightarrow 0.
\tag {4.5}
$$
Since $\overline D\times\overline B\times\overline B'$ is compact, choosing a subsequence of $\{(x_n,y_n,y_n')\}$ if necessary,
one may assume that $\{(x_n,y_n,y_n')\}$ converges to a point $(x_*,y_*, y_*')\in\overline D\times\overline B\times\overline B'$.
Since $y_*,y_*'\not=x_*$ by assumption,
it follows from (4.5) that
$$\displaystyle
\frac{y_*-x_*}{\vert y_*-x_*\vert}+\frac{y_*'-x_*}{\vert y_*'-x_*\vert}=0.
$$
This gives $x_*\in\{sy_*+(1-s)y_*'\,\vert\,0<s<1\}$ and thus $[\overline B\cup
\overline B']\cap\overline D\not=\emptyset$.
This is a contradiction.

Thus, applying (4.4) to the right-hand side of (4.3), we obtain, for a sufficiently large $\tau_0>0$
$$\displaystyle
\inf_{(x,y,y')\in D\times B\times B'}K_{\tau}(x,y,y')
\ge\frac{\tau^2}{2}C_D(B,B'),\,\,\forall\tau\ge\tau_0'.
$$
A combination of this and (4.2) ensures the validity of Lemma 4.1
with
$$\displaystyle
C=\left(\frac{1}{4\pi}\right)^2\frac{\tau^2 C_D(B,B')}{2\text{dist}\,(D,B)\,\text{dist}\,(D,B')}.
$$

\noindent
$\Box$

\proclaim{\noindent Lemma 4.2. }
Assume that $D$ is admissible and its boundary is $C^3$.
Then, there exist positive constants $C$ and $\tau_0$ independent of $\tau$ such that, for all $\tau\ge\tau_0$,
$$\displaystyle
\left\vert\int_{\Bbb R^3\setminus\overline D}
(\nabla\epsilon_f^0\cdot\nabla\epsilon_g^0+\tau^2
\epsilon_f^0\epsilon_g^0)dx-J(\tau;f,g)\right\vert
\le C\tau^{3/2}\int_Ddx\int_{B\times B'} e^{-\tau\phi(x;y,y')}dydy'.
$$

\endproclaim

{\it\noindent Proof.}
Let $0<\delta<\delta_0$.
Let $\phi=\phi_{\delta}$ be a smooth cut-off
function, $0\le\phi(x)\le 1$, and such that: $\phi(x)=1$ if
$d_{\partial D}(x)<\delta$ and $\phi(x)=0$ if $d_{\partial
D}(x)>2\delta$;
$\vert \nabla\phi(x)\vert\le C\delta^{-1}$;
$\vert\nabla^2\phi(x)\vert\le C\delta^{-2}$.

For $*=f,g$
define $\displaystyle (\epsilon_{*}^0)^r(x)=\phi(x)\epsilon_*^0(x^r)$ for $x\in D$ and $v_*^r(x)=\phi(x)v_*(x^r)$ for
$x\in\,\Bbb R^3\setminus\overline D$.

The Lax-Phillips reflection argument starts with the following
expression:
$$\displaystyle
\int_{\Bbb R^3\setminus\overline D}
(\nabla\epsilon_f^0\cdot\nabla\epsilon_g^0+\tau^2\epsilon_f^0\epsilon_g^0)
dx
=J(\tau;f,g)
+\int_{\Bbb R^3\setminus\overline D}\epsilon_f^0(\triangle-\tau^2)v_g^rdx.
\tag {4.6}
$$
In the proof the following relationship between $v_*^r$ and $v_*$
and the boundary condition for $\epsilon_*^0$ are essential:
$$\begin{array}{c}
\displaystyle
\frac{\partial v_*^r}{\partial\nu}
=-\frac{\partial v_*}{\partial\nu}\,\,\text{on}\,\partial D,\\
\\
\displaystyle
v_*^r=v_*=-\epsilon_*^0\,\,\text{on}\,\partial D.
\end{array}
$$

Another device is the following differential identity which is a consequence of (4.15) in
\cite{LP}(see also Appendix 1 in \cite{IE4}):
$$\displaystyle
(\triangle-\tau^2)(v_g^r)
=\phi(x)\sum_{i,j}a_{ij}(x)(\partial_i\partial_jv_g)(x^r)+\sum_{j}b_j(x)(\partial_jv_g)(x^r)+(\triangle\phi)(x)v_g(x^r).
$$
where $a_{ij}(x)$, $i,j=1,2,3$ are $C^1$ in a neighbourhood of $\partial D$, independent of $\phi$ and $v_g$ and satisfy
$$\displaystyle
\exists C>0\,\forall x\in\Bbb R^3\,\,\,\phi(x)\vert a_{ij}(x)\vert\le Cd_{\partial D}(x);
\tag {4.7}
$$
each $b_j(x)$ has the form
$$\displaystyle
b_j(x)=\sum_{jk}b_{jk}(x)\partial_k\phi(x)+d_j(x)\phi(x)
\tag {4.8}
$$
with $b_{jk}(x)$ and $d_j(x)$ which are $C^1$ and $C^0$ in a neighbourhood of $\partial D$, respectively
and independent of $\phi$ and $v_g$.

This together with the change of variables $x=y^r$ yields
$$\begin{array}{c}
\displaystyle
\int_{\Bbb R^3\setminus\overline D}\epsilon_f^0(x)(\triangle-\tau^2)v_g^r(x)dx\\
\\
\displaystyle
=\int_{\Bbb R^3\setminus\overline D}\epsilon_f^0(x)\left\{\phi(x)\sum_{i,j}a_{ij}(x)(\partial_i\partial_jv_g)(x^r)+\sum_{j}b_j(x)(\partial_jv_g)(x^r)
+(\triangle\phi)(x)v_g(x^r)\right\}dx\\
\\
\displaystyle
=\int_{D}\epsilon_f^0(y^r)\left\{\phi(y^r)\sum_{i,j}a_{ij}(y^r)(\partial_i\partial_jv_g)(y)+\sum_{j}b_j(y^r)(\partial_jv_g)(y)
+(\triangle\phi)(y^r)v_g(y)\right\}
J(y)dy,
\end{array}
\tag {4.9}
$$
where $J(y)$ denotes the Jacobian.

Hereafter, we give an estimation for each term in the right-hand side of (4.9) pointwisely,
without making use of integration by parts further.  This idea is exactly same as the proof of
Lemma 3.3 in \cite{LP}.  This is different from the back-scattering case, see also the proof of Lemma 4.2 in \cite{LP}
and Appendix 1 in \cite{IE4} for the comparison.

Since $D$ is admissible we have, for all $y\in D$ with $d_{\partial D}(y)<\delta'$
and all $\tau>\tau_0$
$$\displaystyle
\vert\epsilon_f^0(y^r)\vert
\le C\int_B e^{-\tau\vert y-x\vert}dx,
\tag {4.10}
$$
where $C$ is independent of $y$ and $\tau$.

It is easy to see that, for all $y\in D$,
$$
\displaystyle
\tau^{-2}\vert\partial_i\partial_jv_g(y)\vert
+\tau^{-1}\vert\partial_iv_g(y)\vert
+
\vert v_g(y)\vert\le C\int_{B'}e^{-\tau\vert y-x'\vert}dx'.
$$
From this, (4.10), (4.7) and (4.8) and the choice of $\phi$ we obtain, for all $y\in D$,
$$\begin{array}{c}
\displaystyle
\left\vert\epsilon_f^0(y^r)\left\{\phi(y^r)\sum_{i,j}a_{ij}(y^r)(\partial_i\partial_jv_g)(y)+\sum_{j}b_j(y^r)(\partial_jv_g)(y)
+(\triangle\phi)(y^r)v_g(y)\right\}\right\vert\\
\\
\displaystyle
\le C(\delta\tau^2+\delta^{-1}\tau+\delta^{-2})
\int_{B\times B'}e^{-\tau\phi(y;x,x')}dxdx'.
\end{array}
\tag {4.11}
$$
Choosing $\delta=\tau^{-1/2}$, we have $\delta\tau^2+\delta^{-1}\tau+\delta^{-2}=O(\tau^{3/2})$ as $\tau\longrightarrow\infty$.
Now from (4.9) and (4.11) we obtain the desired conclusion of Lemma 4.2.

\noindent
$\Box$

Thus everything is reduced to studying the asymptotic behaviour
of $J(\tau;f,g)$ as $\tau\longrightarrow\infty$.
For this purpose we employ the asymptotic formula (3.13).
Note that in Section 3 we made use of the formula to give a lower estimate of $J(\tau;f,g)$.
Here using the formula, we determine its leading term as $\tau\longrightarrow\infty$.

From (3.13) we see that the asymptotic behaviour of the following integral
is the key:
$$\displaystyle
\int_{\partial D}\frac{(p-x)\cdot\nu_x}{\vert x-p\vert^2\vert x-p'\vert}
\left(1+\frac{1}{\tau\vert x-p\vert}\right)
e^{-\tau\phi(x;p,p')}dS_x.
$$

\proclaim{\noindent Proposition 4.2.}
Assume that $\Lambda_{\partial D}(p,p')$ is finite and
$$\displaystyle
\text{det}\,(S_q(E_c(p,p'))-S_q(\partial D))>0\,\,\forall q\in\Lambda_{\partial D}(p,p').
$$
Then, we have
$$\begin{array}{c}
\displaystyle
\lim_{\tau\longrightarrow\infty}
\tau e^{\tau\phi(q;p,p')}
\int_{\partial D}
\frac{(p-x)\cdot\nu_x}{\vert x-p\vert^2\vert x-p'\vert}
\left(1+\frac{1}{\tau\vert x-p\vert}\right)
e^{-\tau\phi(x;p,p')}dS_x\\
\\
\displaystyle
=\sum_{q\in\Lambda_{\partial D}(p,p')}
\frac{\pi}{\vert q-p\vert\vert q-p'\vert}
\frac{1}
{\sqrt{\text{det}\,(S_q(E_c(p,p'))-S_q(\partial D))}}.
\end{array}
$$

\endproclaim

Theorem 1.3 is a direct consequence of this together with (3.13) and
(4.1).

In the following subsection, we describe the proof of Proposition 4.2.

\subsection{Proof of Proposition 4.2.}

We employ the Laplace method and so one has to compute the Hessian of the real phase function $\partial D\ni x\longmapsto\phi(x;p,p')$ at all the points
on $\partial D$ where it takes the minimum value.

Let $q\in\Lambda_{\partial D}(p,p')$.
One can choose a local coordinates system
$\sigma=(\sigma_1,\sigma_2)$ around $q$ on $\partial D$
in such a way that $x\in\partial D$ around $q$ has the form
$$\displaystyle
x=x_q(\sigma)=q+\sigma_1\mbox{\boldmath $e$}_1+\sigma_2\mbox{\boldmath $e$}_2+f(\sigma)\nu_q,\,\vert\sigma\vert<\delta,
$$
where $\delta$ is a sufficiently small positive number independent of $q$;
$\mbox{\boldmath $e$}_1$ and $\mbox{\boldmath $e$}_2$ are two unit tangent vectors at $q$ to $\partial D$
which are perpendicular to each other and $\mbox{\boldmath $e$}_1\times\mbox{\boldmath $e$}_2=\nu_q$;
$f=f_q\in C^{3}_0(\Bbb R^3)$ and satisfies $f(0)=0$, $\nabla f(0)=0$; $\nu_x$ takes the form
$$\displaystyle
\nu_x=\frac{-f_1\mbox{\boldmath $e$}_1-f_2\mbox{\boldmath $e$}_2+\nu_q}
{\sqrt{1+f_1^2+f_2^2}},
$$
where $f_1=\partial f/\partial\sigma_1$ and $f_2=\partial f/\partial\sigma_2$;
$dS=\sqrt{1+f_1^2+f_2^2}d\sigma$.

We have, for $z=p,p'$ and $j=1,2$,
$$\displaystyle
\frac{\partial}{\partial\sigma_j}\vert x_q(\sigma)-z\vert
=\frac{1}{\vert x_q(\sigma)-z\vert}
\left\{(x_q(\sigma)-z)\cdot\mbox{\boldmath $e$}_j+f_j(\sigma)(x_q(\sigma)-z)\cdot\nu_q\right\}
\tag {4.12}
$$
and thus
$$\begin{array}{c}
\displaystyle
\frac{\partial^2}{\partial\sigma_k\partial\sigma_j}\vert x_q(\sigma)-z\vert\\
\\
\displaystyle
=-\frac{\{(x_q(\sigma)-z)\cdot\mbox{\boldmath $e$}_k+f_k(\sigma)(x_q(\sigma)-z)\cdot\nu_q\}
\{(x_q(\sigma)-z)\cdot\mbox{\boldmath $e$}_j+f_j(\sigma)(x_q(\sigma)-z)\cdot\nu_q\}}
{\vert x_q(\sigma)-z\vert^3}
\\
\\
\displaystyle
+\frac{1}{\vert x_q(\sigma)-z\vert}
\left(\delta_{kj}+\frac{\partial^2 f}{\partial\sigma_k\partial\sigma_j}(\sigma)(x_q(\sigma)-z)\cdot\nu_q
+f_j(\sigma)f_k(\sigma)\right).
\end{array}
$$
This yields
$$\begin{array}{c}
\displaystyle
\frac{\partial^2}{\partial\sigma_k\partial\sigma_j}\phi(x_q(\sigma);p,p')\vert_{\sigma=0}
=\lambda(q;p,p')\delta_{kj}-a_{kj}(q;p,p'),
\end{array}
\tag {4.13}
$$
where
$$\displaystyle
\lambda(q;p,p')=\frac{1}{\vert p-q\vert}+\frac{1}{\vert p'-q\vert},
$$
$$\begin{array}{c}
\displaystyle
a_{kj}(q;p,p')
=\frac{1}{\vert p-q\vert}
\mbox{\boldmath $A$}\cdot\mbox{\boldmath $e$}_k
\mbox{\boldmath $A$}\cdot\mbox{\boldmath $e$}_j
+\frac{1}{\vert p'-q\vert}
\mbox{\boldmath $A$}'\cdot\mbox{\boldmath $e$}_k
\mbox{\boldmath $A$}'\cdot\mbox{\boldmath $e$}_j
-\frac{\partial^2 f}{\partial\sigma_k\partial\sigma_j}(0)
\left(\mbox{\boldmath $A$}\cdot\nu_q+\mbox{\boldmath $A$}'\cdot\nu_q\right)
\end{array}
$$
and
$$\displaystyle
\mbox{\boldmath $A$}=\mbox{\boldmath $A$}_q(p),\,\,
\mbox{\boldmath $A$}'=\mbox{\boldmath $A$}_q(p').
$$

Note that $\lambda(q;p,p')=\lambda(q;p',p)$ and $a_{kj}(q;p,p')=a_{kj}(q;p',p)$.

The following lemma corresponds to Snell's law in geometrical optics.

\proclaim{\noindent Lemma 4.3.}
Let $p$ and $p'$ be arbitrary points in $\Bbb R^3\setminus\partial D$ such that
$[p,p']\cap\partial D=\emptyset$.
Let $q\in\Lambda_{\partial D}(p,p')$. Then, we have: (i) the vector $\mbox{\boldmath $A$}+\mbox{\boldmath $A$}'$ is parallel to $\nu_q$;
(ii) $\mbox{\boldmath $A$}\cdot\nu_q=\mbox{\boldmath $A$}'\cdot\nu_q\not=0$.
In addition, if $p$ and $p'$ are in $\Bbb R^3\setminus\overline D$, then
we have
$$\displaystyle
\nu_q=-\frac{1}{\displaystyle\sqrt{2(1+\mbox{\boldmath $A$}\cdot\mbox{\boldmath $A$}')}}
(\mbox{\boldmath $A$}+\mbox{\boldmath $A$}'),
\tag {4.14}
$$
that is, the unit inward normal to $E_c(p,p')$ at $q$ with $c=\phi(q;p,p')$ (see (A.17) in Appendix) coincides with
the unit outward normal to $\partial D$ at the same point.

\endproclaim

{\it\noindent Proof.}
Since the function $\displaystyle\sigma\mapsto \phi(x_q(\sigma);p,p')$
takes its minimum at $\sigma=0$, it follows from (4.12) that
$$\begin{array}{c}
\displaystyle
\mbox{\boldmath $A$}\cdot\mbox{\boldmath $e$}_1+\mbox{\boldmath $A$}'\cdot\mbox{\boldmath $e$}_1=0\\
\\
\displaystyle
\mbox{\boldmath $A$}\cdot\mbox{\boldmath $e$}_2+\mbox{\boldmath $A$}'\cdot\mbox{\boldmath $e$}_2=0.
\end{array}
\tag {4.15}
$$
Write
$$\begin{array}{c}
\displaystyle
\mbox{\boldmath $A$}
=\alpha\mbox{\boldmath $e$}_1+\beta\mbox{\boldmath $e$}_2+\gamma\nu_q,\\
\\
\displaystyle
\mbox{\boldmath $A$}'
=\alpha'\mbox{\boldmath $e$}_1+\beta'\mbox{\boldmath $e$}_2+\gamma'\nu_q.
\end{array}
\tag {4.16}
$$
(4.15) is equivalent to $\alpha+\alpha'=0$ and $\beta+\beta'=0$.  Then we have $\gamma^2=\gamma'^2$, that is,
$\gamma=\gamma'$ or $\gamma=-\gamma'$.  Assume that $\gamma=-\gamma'$.   Then from (4.16) we have
$\mbox{\boldmath $A$}+\mbox{\boldmath $A$}'=0$, that is,
$$\displaystyle
\frac{q-p}{\vert q-p\vert}+\frac{q-p'}{\vert q-p'\vert}=0.
$$
This means that $q\in\,[p,p']$ and contradicts the condition $[p,p']\cap\partial D=\emptyset$.
Thus $\gamma=\gamma'$ and we have
$$
\displaystyle
\mbox{\boldmath $A$}+\mbox{\boldmath $A$}'=2\gamma\nu_q
\tag {4.17}
$$
and
$$\displaystyle
\gamma=\pm\frac{\sqrt{1+\mbox{\boldmath $A$}\cdot\mbox{\boldmath $A$}'}}{\sqrt{2}}.
$$
Since $\vert\mbox{\boldmath $A$}+\mbox{\boldmath $A$}'\vert^2=2(1+\mbox{\boldmath $A$}\cdot\mbox{\boldmath $A$}')$,
from the argument above we have $1+\mbox{\boldmath $A$}\cdot\mbox{\boldmath $A$}'>0$.
Thus from (4.17), one gets $\mbox{\boldmath $A$}\cdot\nu_q\not=0$ and $\mbox{\boldmath $A$}'\cdot\nu_q\not=0$.

Note that if $p$ and $p'$ are in $\Bbb R^3\setminus\overline D$, $\gamma$ has to be negative.
The reason is the following.  Assume that $\gamma>0$.  Then $-(\mbox{\boldmath $A$}+\mbox{\boldmath $A$}')$
is directed to $-\nu_q$.
Since $\partial D$ is $C^2$,
one can find a sufficiently small $s>0$ such that $q'\equiv q-s(\mbox{\boldmath $A$}+\mbox{\boldmath $A$}')\in D$.
Since $p\in\Bbb R^3\setminus\overline D$,
one can find a point $p''\in\partial D$ on the segment with endpoints $q'$ and $p$.
Note that $-(\mbox{\boldmath $A$}+\mbox{\boldmath $A$}')$
is directed to the unit inward normal to $E_c(p,p')$.  This together with the condition $c>\vert p-p'\vert$
ensures that both $q'$ and $p$ are in the domain enclosed by $E_c(p,p')$ and thus by the convexity
of the domain, we have $\phi(p'';p,p')<c$.
However, since $p''\in\partial D$, we have $\phi(p'';p,p')\ge c$.  This is a contradiction.
Therefore one gets $\gamma<0$ and now it is easy to see that all the conclusions are valid.

\noindent
$\Box$

It follows from (i) in Lemma 4.3 that
$$\displaystyle
\mbox{\boldmath $A$}\cdot\mbox{\boldmath $e$}_k
\mbox{\boldmath $A$}\cdot\mbox{\boldmath $e$}_j
=\mbox{\boldmath $A$}'\cdot\mbox{\boldmath $e$}_k
\mbox{\boldmath $A$}'\cdot\mbox{\boldmath $e$}_j.
$$
Thus, we have
$$\begin{array}{c}
\displaystyle
\frac{1}{\vert p-q\vert}
\mbox{\boldmath $A$}\cdot\mbox{\boldmath $e$}_k\mbox{\boldmath $A$}\cdot\mbox{\boldmath $e$}_j\\
\\
\displaystyle
=\left(\frac{1}{\vert p-q\vert}+\frac{1}{\vert p'-q\vert}\right)
\mbox{\boldmath $A$}\cdot\mbox{\boldmath $e$}_k
\mbox{\boldmath $A$}\cdot\mbox{\boldmath $e$}_j
-\frac{1}{\vert p'-q\vert}
\mbox{\boldmath $A$}\cdot\mbox{\boldmath $e$}_k\mbox{\boldmath $A$}\cdot\mbox{\boldmath $e$}_j\\
\\
\displaystyle
=\lambda(q;p,p')\mbox{\boldmath $A$}\cdot\mbox{\boldmath $e$}_k
\mbox{\boldmath $A$}\cdot\mbox{\boldmath $e$}_j
-\frac{1}{\vert p'-q\vert}
\mbox{\boldmath $A$}'\cdot\mbox{\boldmath $e$}_k
\mbox{\boldmath $A$}'\cdot\mbox{\boldmath $e$}_j
\end{array}
$$
and hence
$$\displaystyle
\frac{1}{\vert p-q\vert}
\mbox{\boldmath $A$}\cdot\mbox{\boldmath $e$}_k\mbox{\boldmath $A$}\cdot\mbox{\boldmath $e$}_j
+\frac{1}{\vert p'-q\vert}
\mbox{\boldmath $A$}'\cdot\mbox{\boldmath $e$}_k\mbox{\boldmath $A$}'\cdot\mbox{\boldmath $e$}_j
=\lambda(q;p,p')
\mbox{\boldmath $A$}\cdot\mbox{\boldmath $e$}_k
\mbox{\boldmath $A$}\cdot\mbox{\boldmath $e$}_j.
$$
Changing the role of $p$ and $p'$, we also have
$$\displaystyle
\frac{1}{\vert p-q\vert}
\mbox{\boldmath $A$}\cdot\mbox{\boldmath $e$}_k\mbox{\boldmath $A$}\cdot\mbox{\boldmath $e$}_j
+\frac{1}{\vert p'-q\vert}
\mbox{\boldmath $A$}'\cdot\mbox{\boldmath $e$}_k\mbox{\boldmath $A$}'\cdot\mbox{\boldmath $e$}_j
=\lambda(q;p,p')
\mbox{\boldmath $A$}'\cdot\mbox{\boldmath $e$}_k
\mbox{\boldmath $A$}'\cdot\mbox{\boldmath $e$}_j.
$$

From these and (4.14), we obtain another expression for $a_{jk}(q:p,p')$:
$$\begin{array}{c}
\displaystyle
a_{jk}(q;p,p')
=\frac{\lambda(q;p,p')}{2}
\left(\mbox{\boldmath $A$}\cdot\mbox{\boldmath $e$}_k
\mbox{\boldmath $A$}\cdot\mbox{\boldmath $e$}_j
+
\mbox{\boldmath $A$}'\cdot\mbox{\boldmath $e$}_k
\mbox{\boldmath $A$}'\cdot\mbox{\boldmath $e$}_j\right)
+\sqrt{2(1+\mbox{\boldmath $A$}\cdot\mbox{\boldmath $A$}')}\frac{\partial^2 f}{\partial\sigma_k\partial\sigma_j}(0).
\end{array}
$$

This together with (4.13) implies that
$$\begin{array}{c}
\displaystyle
\frac{\partial^2}{\partial\sigma_k\partial\sigma_j}\phi(x_q(\sigma);p,p')\vert_{\sigma=0}
=\sqrt{2(1+\mbox{\boldmath $A$}\cdot\mbox{\boldmath $A$}')}\\
\\
\displaystyle
\times
\left\{
\frac{\lambda(q;p,p')}{\sqrt{2(1+\mbox{\boldmath $A$}\cdot\mbox{\boldmath $A$}')}}
\left(\delta_{jk}
-\frac{1}{2}
(\mbox{\boldmath $A$}\cdot\mbox{\boldmath $e$}_k\mbox{\boldmath $A$}\cdot\mbox{\boldmath $e$}_j
+\mbox{\boldmath $A$}'\cdot\mbox{\boldmath $e$}_k\mbox{\boldmath $A$}'\cdot\mbox{\boldmath $e$}_j)
\right)
-\frac{\partial^2 f}{\partial\sigma_k\sigma_j}(0)\right\}.
\end{array}
\tag {4.18}
$$
From (A.2) we have
$$\begin{array}{c}
\displaystyle
\left(\begin{array}{c} \mbox{\boldmath $e$}_1^T\\
\\
\displaystyle
\mbox{\boldmath $e$}_2^T
\end{array}
\right)S_q(E_c(p,p'))
\left(\begin{array}{cc} \mbox{\boldmath $e$}_1 &
\mbox{\boldmath $e$}_2
\end{array}
\right)\\
\\
\displaystyle
=\frac{\lambda(q;p,p')}{\sqrt{2(1+\mbox{\boldmath $A$}\cdot\mbox{\boldmath $A$}')}}
\left(\delta_{jk}
-\frac{1}{2}
(\mbox{\boldmath $A$}\cdot\mbox{\boldmath $e$}_k\mbox{\boldmath $A$}\cdot\mbox{\boldmath $e$}_j
+\mbox{\boldmath $A$}'\cdot\mbox{\boldmath $e$}_k\mbox{\boldmath $A$}'\cdot\mbox{\boldmath $e$}_j)
\right)_{j\downarrow 1,2;k\rightarrow 1,2}
\end{array}
$$
and we know
$$
\displaystyle
\left(\begin{array}{c} \mbox{\boldmath $e$}_1^T\\
\\
\displaystyle
\mbox{\boldmath $e$}_2^T
\end{array}
\right)S_q(\partial D)
\left(\begin{array}{cc} \mbox{\boldmath $e$}_1 &
\mbox{\boldmath $e$}_2
\end{array}
\right)
=\left(\frac{\partial^2 f}{\partial\sigma_j\partial\sigma_k}(0)\right)_{j\downarrow 1,2;k\rightarrow 1,2}.
\tag {4.19}
$$
Thus we obtain the following formula which gives the geometrical
meaning of the Hessian of $\phi(x_q(\sigma);p,p')$ at $\sigma=0$.

\proclaim{\noindent Lemma 4.4.}
Let $p$ and $p'$ be in $\Bbb R^3\setminus\overline D$ such that
$[p,p']\cap\partial D=\emptyset$.
Let $q\in\Lambda_{\partial D}(p,p')$ and $c=\phi(q;p,p')$.
Then, $c>\vert p-p'\vert$ and we have
$$\displaystyle
\nabla^2_{\sigma}\phi(x_q(\sigma);p,p')\vert_{\sigma=0}
=\sqrt{2(1+\mbox{\boldmath $A$}\cdot\mbox{\boldmath $A$}')}
\left(\begin{array}{c} \mbox{\boldmath $e$}_1^T\\
\\
\displaystyle
\mbox{\boldmath $e$}_2^T
\end{array}
\right)
\left(S_q(E_c(p,p'))-S_q(\partial D)\right)
\left(\begin{array}{cc} \mbox{\boldmath $e$}_1 &
\mbox{\boldmath $e$}_2
\end{array}
\right).
\tag {4.20}
$$

\endproclaim

Since $\phi(x_q(\sigma);p,p')$ takes the minimum at $\sigma=0$, from (4.20) one concludes
that, for all tangent vectors $\mbox{\boldmath $v$}$ at $q$
$$\displaystyle
(S_q(E_c(p,p'))-S_q(\partial D))\mbox{\boldmath $v$}\cdot\mbox{\boldmath $v$}\ge 0.
\tag {4.21}
$$
Thus if (1.10) is satisfied, then from (4.21) one knows that
$S_q(E_c(p,p'))-S_q(\partial D)$ is positive definite on the
tangent space at $q$ and from (4.20)
$$\displaystyle
\text{det}\,\nabla^2_\sigma\phi(x_q(\sigma);p,p')\vert_{\sigma=0}
=2(1+\mbox{\boldmath $A$}\cdot\mbox{\boldmath $A$}')\text{det}\,(S_q(E_c(p,p'))-S_q(\partial D))>0.
$$
And also from (4.14) we have
$$\displaystyle
(p-q)\cdot\nu_q
=\frac{\vert q-p\vert}
{\sqrt{2}}\sqrt{1+\mbox{\boldmath $A$}\cdot\mbox{\boldmath $A$}'}.
$$

Now Proposition 4.2 is a direct consequence of the Laplace method \cite{BH}.

\section{Proof of Theorems 1.2 and 1.4}

\subsection{Extracting the first reflector: proof of Theorem 1.2}

Let $c>\vert p-p'\vert$.
Theorem 1.2 is based on the following proposition which gives a characterization
of a first reflection point $q$ between $p$ and $p'$ in terms of the minimum length
of the broken path connecting $p$ to $q$ and $q$ to $p'+s(q-p')/\vert q-p'\vert$
with a fixed small $s>0$.

\proclaim{\noindent Proposition 5.1.}  Fix $0<s<\eta'$.
Let $c=\min_{x\in\partial D}\phi(x;p,p')$.  We have:

(i)  if $p'+s(\omega;p,p',c)\omega$ belongs to $\partial D$, then
$$\displaystyle
\min_{x\in\partial D}\phi(x;p,p'+s\omega)=c-s;
$$

(ii) if $p'+s(\omega;p,p',c)\omega$ does not belong to $\partial D$, then
$$\displaystyle
\min_{x\in\partial D}\phi(x;p,p'+s\omega)>c-s.
$$
Thus, one has the following characterization of the first reflector:
$$\displaystyle
\Lambda_{\partial D}(p,p')
=\{p'+s(\omega;p,p',c)\omega\,\vert\,\min_{x\in\partial D}\phi(x;p,p'+s\omega)=c-s,\,\,\omega\in S^2\}.
$$

\endproclaim

{\it\noindent Proof.}
Set $\displaystyle p''(\omega)
=p'+s\omega$.
Let $x\in\partial D$.  We have
$$\begin{array}{c}
\displaystyle
\phi(x;p,p''(\omega))
=\vert p-x\vert+\vert x-p''(\omega)\vert\\
\\
\displaystyle
=\vert p-x\vert+\vert (x-p')+(p'-p''(\omega))\vert\\
\\
\displaystyle
\ge\vert p-x\vert+\vert x-p'\vert-\vert p'-p''(\omega)\vert
=\phi(x;p,p')-s.
\end{array}
\tag {5.1}
$$
This gives
$$\displaystyle
\min_{x\in\partial D}\phi(x;p,p''(\omega))\ge c-s.
\tag {5.2}
$$
Now we describe the proof of (i).
Noting that $s<s(\omega;p,p',c)$ and $p'+s(\omega;p,p',c)\omega\in E_c(p,p')$, we have
$$\begin{array}{c}
\displaystyle
\phi(p'+s(\omega;p,p',c)\omega;p,p''(\omega))\\
\\
\displaystyle
=\vert p-(p'+s(\omega;p,p',c)\omega)\vert+\vert (p'+s(\omega;p,p',c)\omega)-p''(\omega)\vert\\
\\
\displaystyle
=\vert p-(p'+s(\omega;p,p',c)\omega)\vert+(s(\omega;p,p',c)-s)\\
\\
\displaystyle
=\vert p-(p'+s(\omega;p,p',c)\omega)\vert+\vert (p'+s(\omega;p,p',c)\omega)-p'\vert-s\\
\\
\displaystyle
=\phi(p'+s(\omega;p,p',c)\omega);p,p')-s
=c-s.
\end{array}
$$
Thus if $p'+s(\omega;p,p',c)\omega$ belongs to $\partial D$, then the inequality in (5.2)
$$\displaystyle
\min_{x\in\partial D}\phi(x;p,p''(\omega))>c-s
$$
never occurs.  Thus it must hold that $\min_{x\in\partial D}\phi(x;p,p''(\omega))=c-s$.

Since we have always (5.2), the statement of (ii) is equivalent to the one that:
if $\min_{x\in\partial D}\phi(x;p,p''(\omega))=c-s$,
then $p'+s(\omega;p,p',c)\omega\in\partial D$.

Assume that $\min_{x\in\partial D}\phi(x;p,p''(\omega))=c-s$.
Choose a point $x\in\partial D$ such that $\phi(x;p,p''(\omega))=c-s$.
Then, from (5.1) we have $\displaystyle
c\ge\phi(x;p,p')$.  Since $c=\min_{x\in\partial D}\phi(x;p,p')$, from this we obtain $\phi(x;p,p')=c$.
Thus we have $\phi(x;p,p''(\omega))=\phi(x;p,p')-s$.  This is equivalent to
$$\displaystyle
\vert x-p''(\omega)\vert=\vert x-p'\vert-s.
\tag {5.3}
$$
Since $x$ is outside the open ball $B''$ centered at $p'$ with
radius $s$, one can find the unique point $x'$ on $\partial B''$
such that $\vert x-x'\vert=\min_{y\in\partial B''}\vert y-x\vert$.
Since we have $\vert x-p'\vert=\vert x-x'\vert+s$, from (5.3) we
obtain $\vert x-p''(\omega)\vert=\vert x-x'\vert$.  Since
$p''(\omega)\in\partial B''$, it must hold that $p''(\omega)=x'$
and thus $x=p'+s(\omega;p,p',c)\omega$. This gives
$p'+s(\omega;p,p',c)\omega\in\partial D$.

\noindent
$\Box$

Now we are ready to describe the proof of Theorem 1.2.

In what follows we denote the open ball centered at a point $z$ and with radius $\rho$ by $B_{\rho}(z)$.
Since $B_{\eta'-s}(p'+s\omega)$ is contained in $B'$,
from $u_f$ on $B'\times]0,\,T[$ with $f=\chi_B$, one gets
$u_f$ on $B_{\eta'-s}(p'+s\omega)\times\,]0,\,T[$.
By Theorem 1.1 we obtain $\inf_{x\in\partial D}\phi(x;p,p'+s\omega)$
for each $\omega\in S^2$.
Thus, from Proposition 5.1 we obtain $\Lambda_{\partial D}(p,p')$ itself.
From formula (4.14) one gets $\nu_q$ at given $q\in\Lambda_{\partial D}(p,p')$.

Thus one can complectly determine the first reflectors between $p$
and $p'$ using the bistatic data $u_f$ on $B'\times\,]0,\,T[$ for
$f=\chi_B$ and sufficiently large and {\it fixed} $T$. In
particular, note that $B$ is {\it fixed}.

{\bf\noindent Remark 5.1.}
We summarize how to detect the points in $\Lambda_{\partial D}(p,p')\subset E_c(p,p')$.

{\bf Step 1.}  Compute $w_f$ on $B'$ with $f=\chi_B$ from the data $u_f$ on $B'\times\,]0,\,T[$.

{\bf Step 2.}  Fix $s$ with $0<s<\eta'$.

{\bf Step 3.}  Choose a direction $\omega\in S^2$.

{\bf Step 4.}  Compute the following integral with
$g=\chi_{B_{\eta'-s}(p'+s\omega)}$:
$$\displaystyle
\int_Bv_gdx-\int_{B_{\eta'-s}(p'+s\omega)}w_fdx.
$$

{\bf Step 5.}
Compute the following quantity:
$$\displaystyle
-\lim_{\tau\longrightarrow\infty}
\frac{1}{\tau}
\log
\left(\int_Bv_gdx-\int_{B_{\eta'-s}(p'+s\omega)}w_fdx\right)-\eta-\eta'+s.
$$

{\bf Step 6.}  If the computed quantity in Step 5 is equal to $c-s$, then
$p'+s(\omega;p,p',c)\omega\in\Lambda_{\partial D}(p,p')$.
If not so, then choose another $\omega$ and go to Step 4.

\subsection{Extracting the geometry of $\partial D$ at a first reflection point: proof of Theorem 1.4}

In this subsection, we consider the case when a point $q\in\Lambda_{\partial D}(p,p')$ is {\it known}.
We use the notation
$\mbox{\boldmath $A$}$ and $\mbox{\boldmath $A$}'$ instead of $\mbox{\boldmath $A$}_q(p)$
and $\mbox{\boldmath $A$}_q(p')$, respectively for simplicity of description.
The aim of this subsection is to extract the geometry of $\partial D$ at $q$.
The main idea is to replace $p'$ in $\Lambda_{\partial D}(p,p')$ with $p'+s\mbox{\boldmath $A$}'$
with a small $s>0$.  The advantage of this idea is described in the following proposition.

\proclaim{\noindent Proposition 5.2.}  Let $c=\min_{x\in\partial D}\phi(x;p,p')$
and satisfy $c>\vert p-p'\vert$.
Fix $0<s<\eta'$.  If $q\in\Lambda_{\partial D}(p,p')$, then
$\displaystyle
\Lambda_{\partial D}\left(p,p'+s\mbox{\boldmath $A$}'\right)=\{q\}$.
Moreover, if $s<\frac{1}{2}(c-\vert p-p'\vert)$, then the operator
$\displaystyle
S_q\left(E_{c-s}\left(p,p'+s\mbox{\boldmath $A$}'\right)\right)
-S_q(\partial D)$
is positive definite on the common tangent space
$T_q(\partial D)=T_q(E_{c-s}(p,p'+s\mbox{\boldmath $A$}'))=T_q(E_c(p,p'))$ at $q$.

\endproclaim

{\it\noindent Proof.}
From the definition of $\Lambda_{\partial D}(p,p')$, $q\in\partial D$ and one can write
$$\displaystyle
q=p'+s(\mbox{\boldmath $A$}';p,p',c)\mbox{\boldmath $A$}'.
$$
Then, we have $\displaystyle\phi(q;p,p'')=c-s$, where
$\displaystyle
p''=p'+s\mbox{\boldmath $A$}'$.  Since $\phi(x;p,p'')\ge c-s$ for all $x\in\partial D$,
this means $c-s=\min_{x\in\partial D}\phi(x;p,p'')$ and thus
$$\displaystyle
q\in\Lambda_{\partial D}\left(p,p'+s\mbox{\boldmath $A$}'\right).
$$
Next let
$$\displaystyle
x\in\Lambda_{\partial D}\left(p,p'+s\mbox{\boldmath $A$}'\right).
$$
We have $\phi(x;p,p'')=\min_{x\in\partial D}\phi(x;p,p'')=c-s$.
Then we have
$$\begin{array}{c}
\displaystyle
c-s=\phi(x;p,p'')
=\vert p-x\vert+\vert x-p''\vert\\
\\
\displaystyle
=\vert p-x\vert+\vert (x-p')+(p'-p'')\vert\\
\\
\displaystyle
\ge\vert p-x\vert+\vert x-p'\vert-\vert p'-p''\vert
=\phi(x;p,p')-s.
\end{array}
$$
Thus $\phi(x;p,p')\le c$ and hence $c=\phi(x;p,p')$. From these we
have $\phi(x;p,p')-s=\phi(x;p,p'')$.  This is equivalent to
$$\displaystyle
\vert x-p''\vert=\vert x-p'\vert-s.
\tag {5.4}
$$
Since $x$ is outside the open ball $B''$ centered at $p'$ with
radius $s$, one can find the unique point $x'$ on $\partial B''$
such that $\vert x-x'\vert=\min_{y\in\partial B''}\vert y-x\vert$.
Since we have $\vert x-p'\vert=\vert x-x'\vert+s$, from (5.4) we
obtain $\vert x-p''\vert=\vert x-x'\vert$.  Since $p''\in\partial
B''$, it must hold that $p''=x'$ and thus $x=p'+s(\mbox{\boldmath
$A$}';p,p',c)\mbox{\boldmath $A$}'=q$.

The last statement is based on the fact that the eigenvectors for both shape operators are common
and $\lambda'>\lambda$, where
$$
\displaystyle
\lambda=\lambda(q;p,p')
=\frac{1}{\vert q-p\vert}+\frac{1}{\vert q-p'\vert},\,\,
\lambda'=\lambda(q;p,p'+s\mbox{\boldmath $A$}')
=\frac{1}{\vert q-p\vert}+\frac{1}{\vert q-p'\vert-s}.
$$
See Appendix for these.  Thus one concludes that the operator
$\displaystyle S_q\left(E_{c-s}\left(p,p'+s\mbox{\boldmath
$A$}'\right)\right) -S_q(E_{c}(p,p'))$ is positive definite.
Since $S_q(E_{c}(p,p'))-S_q(\partial D)\ge 0$, from these one gets
the desired conclusion. Note that the condition
$s<\frac{1}{2}(c-\vert p-p'\vert)$ is just for ensuring that
$\vert p-(p'+s\mbox{\boldmath $A$}')\vert<c-s$.

\noindent
$\Box$

We have
$$\displaystyle
\min_{x\in\partial D,\,y\in\partial B,\,y'\in\partial B_{\eta'-s}(p'+s\mbox{\boldmath $A$}')}\phi(x;y,y')
=\min_{x\in\partial D}\phi(x;p,p'+s\mbox{\boldmath $A$}')-(\eta+\eta'-s)
$$
and, by (i) of Proposition 5.1, $\min_{x\in\partial D}\phi(x;p,p'+s\mbox{\boldmath $A$}')=c-s$ if $q\in\Lambda_{\partial D}(p,p')$.

Thus the condition
$$\displaystyle
T>\min_{x\in\partial D,\,y\in\partial B,\,y'\in\partial B_{\eta'-s}(p'+s\mbox{\boldmath $A$}')}\phi(x;y,y')
$$
is equivalent to (1.5).
Therefore the following proposition is a direct consequence of Theorem 1.3 together
with Proposition 5.2.

\proclaim{\noindent Proposition 5.3.}
Fix $0<s<\min\,(\eta',(c-\vert p-p'\vert)/2)$.
Let $B$ and $B'$ satisfy $[\overline B\cup\overline B']\cap\overline D=\emptyset$.
Let $f=\chi_B$ and $g=\chi_{B_{\eta'-s}(p'+s\mbox{\boldmath $A$}')}$.
Let $T$ satisfy (1.5).
Let $q\in\Lambda_{\partial D}(p,p')$ and set $c=\phi(q;p,p')$.

If $D$ is admissible and $\partial D$ is $C^3$, then we have
$$\begin{array}{c}
\displaystyle
\lim_{\tau\longrightarrow\infty}
\tau^4e^{\tau(c-\eta-\eta')}
\left(\int_{B}v_gdx-\int_{B_{\eta'-s}(p'+s\mbox{\boldmath $A$}')}w_fdx\right)\\
\\
\displaystyle
=\frac{\pi}{2}\left(\frac{\eta}{\vert q-p\vert}\right)\cdot
\left(\frac{\eta'-s}{\vert q-p'\vert-s}\right)\cdot
\frac{1}
{\sqrt{\text{det}\,(S_q(E_{c-s}(p,p'+s\mbox{\boldmath $A$}'))-S_q(\partial D))}}.
\end{array}
\tag {5.5}
$$

\endproclaim

It is quite interesting to know the quantities contained in
$\text{det}\,(S_q(E_{c-s}(p,p'+s\mbox{\boldmath
$A$}'))-S_q(\partial D))$. The following lemma whose proof is
given in Appendix clarifies them.

\proclaim{\noindent Lemma 5.1.}
Let $0\le s<\min\,(\eta',(c-\vert p-p'\vert)/2)$.
Let $q\in\Lambda_{\partial D}(p,p')$ and set $c=\phi(q;p,p')$.
One has
$$\begin{array}{c}
\displaystyle
\text{det}\,(S_q(E_{c-s}(p,p'+s\mbox{\boldmath $A$}'))-S_q(\partial D))
=
\frac{\lambda(q;p,p'+s\mbox{\boldmath $A$}')^2}{4}\\
\\
\displaystyle
-\sqrt{\frac{2}{1+\mbox{\boldmath $A$}\cdot\mbox{\boldmath $A$}'}}\lambda(q;p,p'+s\mbox{\boldmath $A$}')
\left(H_{\partial D}(q)-\frac{S_{q}(\partial D)(\mbox{\boldmath $A$}\times\mbox{\boldmath $A$}')
\cdot
(\mbox{\boldmath $A$}\times\mbox{\boldmath $A$}')}{2(1+\mbox{\boldmath $A$}\cdot\mbox{\boldmath $A$}')}\right)
+K_{\partial D}(q).
\end{array}
\tag {5.6}
$$
\endproclaim

Now we are ready to describe the proof of Theorem 1.4.

Choose $0<s_1<s_2<\min\,(\eta',(c-\vert p-p'\vert)/2)$.
Let $s=s_1, s_2$.
By (5.5) in Proposition 5.3, one gets $\text{det}\,(S_q(E_{c-s}(p,p'))-S_q(\partial D))$ for $s=s_1,s_2$.
From (5.6) we obtain the system
$$\left(\begin{array}{lr}
\displaystyle
-\sqrt{\frac{2}{1+\mbox{\boldmath $A$}\cdot\mbox{\boldmath $A$}'}}\lambda(q;p,p'+s_1\mbox{\boldmath $A$}')
& 1\\
\\
\displaystyle
-\sqrt{\frac{2}{1+\mbox{\boldmath $A$}\cdot\mbox{\boldmath $A$}'}}\lambda(q;p,p'+s_2\mbox{\boldmath $A$}')
& 1
\end{array}
\right)
\mbox{\boldmath $X$}
=\mbox{\boldmath $F$},
\tag {5.7}
$$
where
$$\displaystyle
\mbox{\boldmath $X$}
=\left(\begin{array}{c}
\displaystyle
H_{\partial D}(q)-\frac{S_{q}(\partial D)(\mbox{\boldmath $A$}\times\mbox{\boldmath $A$}')
\cdot
(\mbox{\boldmath $A$}\times\mbox{\boldmath $A$}')}{2(1+\mbox{\boldmath $A$}\cdot\mbox{\boldmath $A$}')}\\
\\
\displaystyle
K_{\partial D}(q)
\end{array}
\right)
$$
and
$$\displaystyle
\mbox{\boldmath $F$}
=\left(\begin{array}{c}
\displaystyle
\text{det}\,(S_q(E_{c-s_1}(p,p'+s_1\mbox{\boldmath $A$}'))-S_q(\partial D))-\frac{\lambda(q;p,p'+s_1\mbox{\boldmath $A$}')^2}{4}
\\
\\
\displaystyle
\text{det}\,(S_q(E_{c-s_2}(p,p'+s_2\mbox{\boldmath $A$}'))-S_q(\partial D))-\frac{\lambda(q;p,p'+s_2\mbox{\boldmath $A$}')^2}{4}
\end{array}
\right).
$$
Since $\lambda(q;p,p'+s_1\mbox{\boldmath $A$}')<\lambda(q;p,p'+s_2\mbox{\boldmath $A$}')$,
(5.7) is uniquely solvable with respect to $\mbox{\boldmath $X$}$.

This completes the proof of Theorem 1.4.

{\bf\noindent Remark 5.2.}
It follows from (5.6) that
$$\begin{array}{c}
\displaystyle
\text{det}\,(S_q(E_c(p,p'))-S_q(\partial D))
-P_{\partial D}(\mu;q)\\
\\
\displaystyle
=-\mu\left\{\left(\sqrt{\frac{2}{1+\mbox{\boldmath $A$}\cdot\mbox{\boldmath $A$}'}}-1\right)
2H_{\partial D}(q)
-\sqrt{\frac{2}{1+\mbox{\boldmath $A$}\cdot\mbox{\boldmath $A$}'}}
\frac{S_{q}(\partial D)(\mbox{\boldmath $A$}\times\mbox{\boldmath $A$}')
\cdot
(\mbox{\boldmath $A$}\times\mbox{\boldmath $A$}')}{(1+\mbox{\boldmath $A$}\cdot\mbox{\boldmath $A$}')}\right\},
\end{array}
$$
where
$$\displaystyle
\mu=\frac{\lambda(q;p,p')}{2}
=\frac{1}{2}\left(\frac{1}{\vert q-p\vert}+\frac{1}{\vert q-p'\vert}\right).
$$
Since this right-hand side has a bound
$O(\vert\mbox{\boldmath $A$}\times\mbox{\boldmath $A$}'\vert^2)$,
this formula indicates an effect of the bistatic data on $\text{det}\,(S_q(E_c(p,p'))-S_q(\partial D))$.

\section{Summary and some of open problems}

This paper is concerned with an inverse obstacle problem
which employs the dynamical scattering data of acoustic wave
over a finite time interval.
The unknown obstacle $D$ is assumed to be {\it sound-soft} one.
The governing equation of the wave is given by the classical wave equation.
The wave is generated by the initial data which is a characteristic function
of an open ball $B$ centered at $p$ and observed
over a finite time interval on a different ball $B'$ centered at $p'$.
It is assumed that $[\overline B\cup\overline{B'}]\cap\overline D=\emptyset$.
The observed data are the so-called {\it bistatic data}.
This is a simple mathematical model of the data collection process using an acoustic wave/electromagnetic wave
such as, bistatic active {\it sonar}, {\it radar}, etc.
This paper aims at developing an enclosure method which employs the
bistatic data.

It is shown that from the data with some additional assumptions on the lower bound
of $T$ one can extract:

(i) the {\it first arrival time} in the
geometrical optics sense, that is, the shortest length of the broken paths
connecting $p$ to a point $q\in\partial D$ and $q$ to $p$;

(ii) the {\it first reflection points} between $p$ and $p'$, that is, all the points $q\in\partial D$ that minimize the length
of the broken paths connecting $p$ to $q$ and $q$ to $p$.

(iii)  the tangent planes of $\partial D$ at all the first reflection points.

It is also shown that, under the {\it admissibility} condition
for $D$, one can extract the Gauss curvature at an arbitrary first
reflection point and the mean curvature with an additional term
which depends on the positions of $p$, $p'$ and the first
reflection point.  As a byproduct, for an example, a constructive proof of a uniqueness theorem for a spherical obstacle
using the bistatic data is also given.

We think that the problem taken up in this paper is a prototype of
other various interesting problems.  It is quite interesting
whether the approach presented here can be applied to them or to develope its necessary modification. Here
we mention some of them.

$\bullet$  Consider the sound-hard obstacle case or the obstacles with a dissipative boundary condition (cf. \cite{PP}).
And also it is quite important to consider the corresponding problem
for the {\it Maxwell system}. These remain open.

$\bullet$  It would be interesting to consider also the case when
obstacles are embedded in one of the {\it two layers} with known
different propagation speeds and both the source and receivers are
placed in another layer.

$\bullet$  Maybe the most interesting problem is that of extracting geometrical information about an unknown
obstacle $D$ {\it behind} a known obstacle $D_0$ from the monostatic or bistatic data over a finite time interval.
$B$ and $B'$ satisfy $[\overline B\cup\overline{B'}]\cap(\overline D_0\cup\overline D)=\emptyset$ and
$D$ is optically invisible from
$B$ and $B'$ because of the existence of $D_0$.

$\bullet$  It would be interesting to consider the case when $B'$ is placed in the {\it shadow region} of $D$ with respect to $B$.
It means that $B'$ is {\it invisible} from $B$ because of the existence of $D$ between them.
This is the case when $[\overline {B}\cup\overline{B'}]\cap\partial D\not=\emptyset$.
What information about $D$ can one extract from $u_f$ on $B'\times\,]0,\,T[$?

Although the author's interest is pursuit of the possibility of the enclosure method itself, 
research by other approaches to the problem taken up in this paper is also expected.
Someone may think about the use of {\it geometrical optics} in the {\it time domain} as Majda has done in \cite{M}
for the problem considered in \cite{LP}.  See also pages 440-447 in \cite{T} for geometrical optics in the time domain
and \cite{PS} for Majda's approach.
His approach heavily depends on the hyperbolic nature of the governing equation
in contrast to our approach and it should be noted that the existing results by our approach can cover inverse problems
for different type of equations like elliptic, parabolic or hyperbolic ones.
Since this paper has not aimed at the comparative study of various approaches,
we leave it to other opportunities.

$$\quad$$

\centerline{{\bf Acknowledgement}}

This research was partially supported by Grant-in-Aid for
Scientific Research (C)(No. 21540162) of Japan  Society for the
Promotion of Science.

\section{Appendix}

\subsection{Proof of (1.7)}

First we prove that
$$\displaystyle
\min_{x\in\partial D}\phi(x;p,p')-(\eta+\eta')
\ge
\min_{x\in\partial D,\,y\in\partial B,\,y'\in\partial B'}\phi(x;y,y').
$$
Choose $q\in\partial D$ such that
$$\displaystyle
\phi(q;p,p')=\min_{x\in\partial D}\phi(x;p,p').
$$
One can find $y_0\in\partial B\cap [q,p]$ and $y_0'\in\partial B'\cap [q,p']$.
Since $\vert y_0-q\vert=\vert p-q\vert-\eta$ and $\vert y_0'-q\vert=\vert p'-q\vert-\eta'$,
we have
$$\displaystyle
\phi(q;y_0,y_0')
=\phi(q;p,p')-(\eta+\eta').
$$
This yields
$$\displaystyle
\phi(q;p,p')-(\eta+\eta')\ge\min_{x\in\partial D,\,y\in\partial B,\,y'\in\partial B'}\phi(x;y,y').
$$

Next we prove that
$$\displaystyle
\min_{x\in\partial D}\phi(x;p,p')-(\eta+\eta')
\le
\min_{x\in\partial D,\,y\in\partial B,\,y'\in\partial B'}\phi(x;y,y').
$$

Choose $q\in\partial D$, $y_0\in\partial B$ and $y_0'\in\partial B'$ such that
$$
\displaystyle
\phi(q;y_0,y_0')
=\min_{x\in\partial D,\,y\in\partial B,\,y'\in\partial B'}
\phi(x;y,y').
$$
Let $y(s)$ be an arbitrary curve on $\partial B$ such that $y(0)=y_0$.
Since $\phi(q;y(s),y_0')$ takes its minimum value at $s=0$, we have
$$\displaystyle
\frac{d}{ds}\phi(q;y(s),y_0')\vert_{s=0}=0.
$$
This gives
$$\displaystyle
\frac{q-y_0}{\vert q-y_0\vert}\cdot\frac{dy}{ds}(0)=0.
$$
Since $dy/ds(0)$ can be an arbitrary vector perpendicular to
the normal vector at $y_0\in\partial B$ and $q$ is outside of $B$,
we have
$$\displaystyle
\frac{q-y_0}{\vert q-y_0\vert}
=\frac{y_0-p}{\vert y_0-p\vert}.
$$
This yields $y_0\in\,[q,p]$ and similarly we have
$y_0'\in\,[q,p']$. Thus we have $\vert y_0-q\vert=\vert
p-q\vert-\eta$ and $\vert y_0'-q\vert=\vert p'-q\vert-\eta'$. This
yields
$$\displaystyle
\phi(q;y_0,y_0')=\phi(q;p,p')-(\eta+\eta')
$$
and thus
$$\displaystyle
\phi(q;y_0,y_0')\ge\min_{x\in\partial D}\phi(x;p,p')-(\eta+\eta').
$$

\subsection{Proof of Proposition 1.1.}

{\it\noindent Proof of (i).}  Since every point on $\partial D$ is a point of support of $D$, we have (3.6).
And from (3.6) we have $\epsilon_f^0(y^r)\ge -v_f(y)$ for $d_{\partial D}(y)<<1$ and $y\in D$.
Since $v_f\ge 0$ and $\epsilon_f^0<0$ (maximum principle),  we have, for $d_{\partial D}(y)<<1$ and
$y\in D$, $\displaystyle
\vert\epsilon_f^0(y^r)\vert\le v_f(y)$. This yields the desired estimate.

Let $c=\min_{x\in\partial D}\phi(x;p,p')$.
Assume that there exist two distinct points $q$ and $q'$ on $\Lambda_{\partial D}(p,p')$.
Since $D$ is convex, $\overline D$ becomes also convex.  Thus every point
on $[q,q']$ are in $\overline D$.  Choose an arbitrary point $q''$ on $[q,q']\setminus\{q,q'\}$.
We have $q''\in\overline D$.
Since both $q$ and $q'$ are on $E_c(p,p')$ with $c=\min_{x\in\partial D}\phi(x;p,p')$,
it is clear that every point $y$ on $[q,q']\setminus\{q,q'\}$ satisfies
$\phi(y;p,p')<c$ and thus $\phi(q'';p,p')<c$.  Since $\phi(x;p,p')\ge c$ for all $x\in\partial D$
and $q''\in D\cup\partial D$,
it must hold that $q''\in D$.
Since $p$ in $\Bbb R^3\setminus\overline D$, there exists a point $q'''$ on $[q'', p]$ such that $q'''\in\partial D$.
Then it is clear to have $\phi(q''';p,p')\le \phi(q'';p,p')$
and thus $\phi(q''';p,p')<c$.  However, since $q'''\in\partial D$, we have
$\phi(q''';p,p')\ge c$.  This is a contradiction.

{\it\noindent Proof of (ii).}
Since $\partial D$ is in the half space $(x-q)\cdot\nu_q\le 0$, it is easy to see that
$S_q(\partial D)\le 0$ as the quadratic form on $T_q(\partial D)$.  On the other hand,
we have $S_q(E_c(p,p'))$ is positive definite as the quadratic form on $T_q(E_c(p,p'))=T_q(\partial D)$
(see Appendix).
Thus $S_q(E_c(p,p'))-S_q(\partial D)$ is positive definite as the quadratic form on the common tangent space
and this yields (1.14).

\subsection{The shape operator for a spheroid}

Since $c>\vert p-p'\vert$, we have $[p,p']\cap E_c(p,p')=\emptyset$.
Thus, for all $x\in E_c(p,p')$,
$$\displaystyle
\frac{x-p}{\vert x-p\vert}+\frac{x-p'}{\vert x-p'\vert}\not=0
$$
and hence
$$\displaystyle
1+\frac{x-p}{\vert x-p\vert}\cdot\frac{x-p'}{\vert x-p'\vert}>0.
$$
Since
$$\displaystyle
\nabla\phi(x;p,p')
=\frac{x-p}{\vert x-p\vert}+\frac{x-p'}{\vert x-p'\vert},
$$
we conclude that $E_c(p,p')$ is a $C^{\infty}$ surface and clearly compact.
Let $\nu_x$ denote the unit {\it inward normal} for $x\in E_c(p,p')$.
We have
$$\displaystyle
\nu_x=-\frac{\nabla\phi}{\vert\nabla\phi\vert}
$$
and note that
$$\displaystyle
\vert\nabla\phi\vert
=\sqrt{2}\sqrt{1+\frac{x-p}{\vert x-p\vert}\cdot\frac{x-p'}{\vert x-p'\vert}}.
$$

These give the expression
$$\displaystyle
-\nu_x=
\frac{\mbox{\boldmath $A$}(x)+\mbox{\boldmath $A$}'(x)}{\sqrt{2(1+\mbox{\boldmath $A$}(x)\cdot\mbox{\boldmath $A$}'(x))}},
\tag {A.1}
$$
where
$$\displaystyle
\mbox{\boldmath $A$}(x)=\frac{x-p}{\vert x-p\vert},\,\,\mbox{\boldmath $A$}'(x)=\frac{x-p'}{\vert x-p'\vert}.
$$

Now we are ready to prove the following proposition.

\proclaim{\noindent Proposition A.1.} Let $S_x$ denote the shape
operator at $x\in E_c(p,p')$ with respect to $\nu_x$ which is the
unit {\it inward normal} to $E_c(p,p')$. We have, for all
$\mbox{\boldmath $v$}\in T_xE_c(p,p')$,
$$\displaystyle
S_x(\mbox{\boldmath $v$})
=\frac{\lambda(x)}{\sqrt{2(1+\mbox{\boldmath $A$}(x)\cdot\mbox{\boldmath $A$}'(x))}}
\left(I_3-\frac{1}{2}\mbox{\boldmath $A$}(x)\otimes\mbox{\boldmath $A$}(x)
-\frac{1}{2}\mbox{\boldmath $A$}'(x)\otimes\mbox{\boldmath $A$}'(x)\right)\mbox{\boldmath $v$},
\tag {A.2}
$$
where
$$\displaystyle
\lambda(x)=\frac{1}{\vert x-p\vert}+\frac{1}{\vert x-p'\vert}.
$$

\endproclaim

{\it\noindent Proof.}
Let $x=x(\sigma_1,\sigma_2)$ be an equation for $E_c(p,p')$ around $q\in E_c(p,p')$.
It means that $x(0,0)=q$ and $\phi(x(\sigma_1,\sigma_2);p,p')=c$.
Since
$$\displaystyle
\frac{\partial}{\partial\sigma_j}\vert x-p\vert
=\frac{x-p}{\vert x-p\vert}\cdot\frac{\partial x}{\partial\sigma_j},
$$
we have
$$\begin{array}{c}
\displaystyle
\frac{\partial}{\partial\sigma_j}\mbox{\boldmath $A$}(x)
=\frac{1}{\vert x-p\vert}\frac{\partial x}{\partial\sigma_j}
-\frac{x-p}{\vert x-p\vert^2}\frac{x-p}{\vert x-p\vert}\cdot\frac{\partial x}{\partial\sigma_j}\\
\\
\displaystyle
=\frac{1}{\vert x-p\vert}\left\{\frac{\partial x}{\partial\sigma_j}
-\mbox{\boldmath $A$}(x)\left(\mbox{\boldmath $A$}(x)\cdot\frac{\partial x}{\partial\sigma_j}\right)\right\}\\
\\
\displaystyle
=\frac{1}{\vert x-p\vert}
\left(I_3-\mbox{\boldmath $A$}(x)\otimes\mbox{\boldmath $A$}(x)\right)\frac{\partial x}{\partial\sigma_j}.
\end{array}
$$
This together with a corresponding expression for
$(\partial/\partial\sigma_j)\mbox{\boldmath $A$}'(x)$ gives
$$\displaystyle
\frac{\partial}{\partial\sigma_j}(\mbox{\boldmath $A$}(x)+\mbox{\boldmath $A$}'(x))
=
\left(\lambda(x)I_3-\frac{1}{\vert x-p\vert}\mbox{\boldmath $A$}(x)\otimes\mbox{\boldmath $A$}(x)-
\frac{1}{\vert x-p'\vert}\mbox{\boldmath $A$}'(x)\otimes\mbox{\boldmath $A$}'(x)\right)\frac{\partial x}{\partial\sigma_j}.
$$
Moreover, since
$$\displaystyle
\frac{\partial x}{\partial\sigma_j}\cdot\mbox{\boldmath $A$}(x)=-\frac{\partial x}{\partial\sigma_j}\cdot\mbox{\boldmath $A$}'(x),
$$
we have
$$\begin{array}{c}
\displaystyle
\frac{\partial}{\partial\sigma_j}(\mbox{\boldmath $A$}(x)\cdot\mbox{\boldmath $A$}'(x))\\
\\
\displaystyle
=\frac{1}{\vert x-p\vert}
\left(I_3-\mbox{\boldmath $A$}(x)\otimes\mbox{\boldmath $A$}(x)\right)\frac{\partial x}{\partial\sigma_j}\cdot\mbox{\boldmath $A$}'(x)
+
\frac{1}{\vert x-p'\vert}
\mbox{\boldmath $A$}(x)\cdot\left(I_3-\mbox{\boldmath $A$}'(x)\otimes\mbox{\boldmath $A$}'(x)\right)\frac{\partial x}{\partial\sigma_j}\\
\\
\displaystyle
=\frac{1}{\vert x-p\vert}\frac{\partial x}{\partial\sigma_j}\cdot\mbox{\boldmath $A$}'(x)
+\frac{1}{\vert x-p'\vert}\frac{\partial x}{\partial\sigma_j}\cdot\mbox{\boldmath $A$}(x)\\
\\
\displaystyle
-\frac{1}{\vert x-p\vert}\left(\mbox{\boldmath $A$}(x)\cdot\frac{\partial x}{\partial\sigma_j}\right)\left(\mbox{\boldmath $A$}(x)
\cdot\mbox{\boldmath $A$}'(x)\right)
-\frac{1}{\vert x-p'\vert}\left(\mbox{\boldmath $A$}'(x)\cdot\frac{\partial x}{\partial\sigma_j}\right)
\left(\mbox{\boldmath $A$}(x)\cdot\mbox{\boldmath $A$}'(x)\right)\\
\\
\displaystyle
=-\frac{1}{\vert x-p\vert}(1+\mbox{\boldmath $A$}(x)\cdot\mbox{\boldmath $A$}'(x))\left(\mbox{\boldmath $A$}(x)
\cdot\frac{\partial x}{\partial\sigma_j}\right)
-\frac{1}{\vert x-p'\vert}(1+\mbox{\boldmath $A$}(x)\cdot\mbox{\boldmath $A$}'(x))
\left(\mbox{\boldmath $A$}'(x)\cdot\frac{\partial x}{\partial\sigma_j}\right)\\
\\
\displaystyle
=-(1+\mbox{\boldmath $A$}(x)\cdot\mbox{\boldmath $A$}'(x))
\left(\frac{1}{\vert x-p\vert}\mbox{\boldmath $A$}(x)\cdot\frac{\partial x}{\partial\sigma_j}
+\frac{1}{\vert x-p'\vert}\mbox{\boldmath $A$}'(x)\cdot\frac{\partial x}{\partial\sigma_j}\right)
\end{array}
$$
we obtain
$$\begin{array}{c}
\displaystyle
\frac{\partial}{\partial\sigma_j}
\frac{1}{\sqrt{1+\mbox{\boldmath $A$}(x)\cdot\mbox{\boldmath $A$}'(x)}}\\
\\
\displaystyle
=-\frac{1}{2}
\frac{\displaystyle\frac{\partial}{\partial\sigma_j}(\mbox{\boldmath $A$}(x)\cdot\mbox{\boldmath $A$}'(x))}
{\displaystyle
(1+\mbox{\boldmath $A$}(x)\cdot\mbox{\boldmath $A$}'(x))^{3/2}}\\
\\
\displaystyle
=\frac{1}{2(1+\mbox{\boldmath $A$}(x)\cdot\mbox{\boldmath $A$}'(x))^{1/2}}
\left(\frac{1}{\vert x-p\vert}\mbox{\boldmath $A$}(x)\cdot\frac{\partial x}{\partial\sigma_j}
+\frac{1}{\vert x-p'\vert}\mbox{\boldmath $A$}'(x)\cdot\frac{\partial x}{\partial\sigma_j}\right).
\end{array}
$$
Thus one gets
$$\begin{array}{c}
\displaystyle
\frac{\partial}{\partial\sigma_j}
\frac{\mbox{\boldmath $A$}(x)+\mbox{\boldmath $A$}'(x)}{\sqrt{1+\mbox{\boldmath $A$}(x)\cdot\mbox{\boldmath $A$}'(x)}}
=\frac{\mbox{\boldmath $A$}(x)+\mbox{\boldmath $A$}'(x)}
{2\sqrt{1+\mbox{\boldmath $A$}(x)\cdot\mbox{\boldmath $A$}'(x)}}\left(\frac{1}{\vert x-p\vert}\frac{\partial x}{\partial\sigma_j}
\cdot\mbox{\boldmath $A$}(x)
+\frac{1}{\vert x-p'\vert}\frac{\partial x}{\partial\sigma_j}\cdot\mbox{\boldmath $A$}'(x)\right)
\\
\\
\displaystyle
+\frac{1}{\sqrt{1+\mbox{\boldmath $A$}(x)\cdot\mbox{\boldmath $A$}'(x)}}
\left(\lambda(x)I_3-\frac{1}{\vert x-p\vert}\mbox{\boldmath $A$}(x)\otimes\mbox{\boldmath $A$}(x)-\frac{1}{\vert x-p'\vert}
\mbox{\boldmath $A$}'(x)\otimes\mbox{\boldmath $A$}'(x)\right)\frac{\partial x}{\partial\sigma_j}
\end{array}
$$
and thus
$$\begin{array}{c}
\displaystyle
\sqrt{1+\mbox{\boldmath $A$}(x)\cdot\mbox{\boldmath $A$}'(x)}\frac{\partial}{\partial\sigma_j}
\frac{\mbox{\boldmath $A$}(x)+\mbox{\boldmath $A$}'(x)}{\sqrt{1+\mbox{\boldmath $A$}(x)\cdot\mbox{\boldmath $A$}'(x)}}\\
\\
\displaystyle
=\frac{1}{2}(\mbox{\boldmath $A$}(x)+\mbox{\boldmath $A$}'(x))\frac{1}{\vert x-p\vert}\frac{\partial x}{\partial\sigma_j}
\cdot\mbox{\boldmath $A$}(x)
+\frac{1}{2}(\mbox{\boldmath $A$}(x)+\mbox{\boldmath $A$}'(x))\frac{1}{\vert x-p'\vert}\frac{\partial x}{\partial\sigma_j}
\cdot\mbox{\boldmath $A$}'(x)\\
\\
\displaystyle
-\mbox{\boldmath $A$}(x)\frac{1}{\vert x-p\vert}\frac{\partial x}{\partial\sigma_j}\cdot\mbox{\boldmath $A$}(x)
-\mbox{\boldmath $A$}'(x)\frac{1}{\vert x-p'\vert}\frac{\partial x}{\partial\sigma_j}\cdot\mbox{\boldmath $A$}'(x)
+\lambda(x)\frac{\partial x}{\partial\sigma_j}\\
\\
\displaystyle
=-\frac{1}{2}\mbox{\boldmath $A$}(x)\frac{1}{\vert x-p\vert}\frac{\partial x}{\partial\sigma_j}\cdot\mbox{\boldmath $A$}(x)
-\frac{1}{2}\mbox{\boldmath $A$}'(x)\frac{1}{\vert x-p'\vert}\frac{\partial x}{\partial\sigma_j}\cdot\mbox{\boldmath $A$}'(x)\\
\\
\displaystyle
+\frac{1}{2}\mbox{\boldmath $A$}'(x)\frac{1}{\vert x-p\vert}\frac{\partial x}{\partial\sigma_j}\cdot\mbox{\boldmath $A$}(x)
+\frac{1}{2}\mbox{\boldmath $A$}(x)\frac{1}{\vert x-p'\vert}\frac{\partial x}{\partial\sigma_j}\cdot\mbox{\boldmath $A$}'(x)
+\lambda(x)\frac{\partial x}{\partial\sigma_j}
\\
\\
\displaystyle
=-\frac{1}{2}\mbox{\boldmath $A$}(x)\frac{1}{\vert x-p\vert}\frac{\partial x}{\partial\sigma_j}\cdot\mbox{\boldmath $A$}(x)
-\frac{1}{2}\mbox{\boldmath $A$}'(x)\frac{1}{\vert x-p'\vert}\frac{\partial x}{\partial\sigma_j}\cdot\mbox{\boldmath $A$}'(x)\\
\\
\displaystyle
-\frac{1}{2}\mbox{\boldmath $A$}'(x)\frac{1}{\vert x-p\vert}\frac{\partial x}{\partial\sigma_j}\cdot\mbox{\boldmath $A$}'(x)
-\frac{1}{2}\mbox{\boldmath $A$}(x)\frac{1}{\vert x-p'\vert}\frac{\partial x}{\partial\sigma_j}\cdot\mbox{\boldmath $A$}(x)
+\lambda(x)\frac{\partial x}{\partial\sigma_j}
\\
\\
\displaystyle
=-\frac{1}{2}\mbox{\boldmath $A$}(x)\lambda(x)\frac{\partial x}{\partial\sigma_j}\cdot\mbox{\boldmath $A$}(x)
-\frac{1}{2}\mbox{\boldmath $A$}'(x)\lambda(x)\frac{\partial x}{\partial\sigma_j}\cdot\mbox{\boldmath $A$}'(x)
+\lambda(x)\frac{\partial x}{\partial\sigma_j}.
\end{array}
$$
Therefore we obtain
$$\begin{array}{c}
\displaystyle
\frac{\partial}{\partial\sigma_j}
\frac{\mbox{\boldmath $A$}(x)+\mbox{\boldmath $A$}'(x)}{\sqrt{1+\mbox{\boldmath $A$}(x)\cdot\mbox{\boldmath $A$}'(x)}}\\
\\
\displaystyle
=\frac{\lambda(x)}{\sqrt{1+\mbox{\boldmath $A$}(x)\cdot\mbox{\boldmath $A$}'(x)}}
\left(I_3-\frac{1}{2}\mbox{\boldmath $A$}(x)\otimes\mbox{\boldmath $A$}(x)-
\frac{1}{2}\mbox{\boldmath $A$}'(x)\otimes\mbox{\boldmath $A$}'(x)\right)\frac{\partial x}{\partial\sigma_j}
\end{array}
$$
and hence
$$\begin{array}{c}
\displaystyle
-\frac{\partial}{\partial\sigma_j}\nu_x
=\frac{\lambda(x)}{\sqrt{2(1+\mbox{\boldmath $A$}(x)\cdot\mbox{\boldmath $A$}'(x))}}
\left(I_3-\frac{1}{2}\mbox{\boldmath $A$}(x)\otimes\mbox{\boldmath $A$}(x)-
\frac{1}{2}\mbox{\boldmath $A$}'(x)\otimes\mbox{\boldmath $A$}'(x)\right)\frac{\partial x}{\partial\sigma_j}.
\end{array}
$$
Since
$$\displaystyle
S_x\left(\xi_1\frac{\partial x}{\partial\sigma_1}\vert_{\sigma=0}+\xi_2\frac{\partial x}{\partial\sigma_2}\vert_{\sigma=0}\right)
=-\left(\xi_1\frac{\partial}{\partial\sigma_1}\nu_x\vert_{\sigma=0}+\xi_2\frac{\partial}{\partial\sigma_2}\nu_x\vert_{\sigma=0}\right),
$$
we obtain (A.2).

\noindent
$\Box$

From (A.2) one can compute the principle curvatures at $x\in E_c(p,p')$.

Consider the case $\mbox{\boldmath $A$}(x)\not=\mbox{\boldmath $A$}'(x)$.  Since $\mbox{\boldmath $A$}(x)$ and
$\mbox{\boldmath $A$}'(x)$ are unit vectors and $\mbox{\boldmath $A$}(x)\not=-\mbox{\boldmath $A$}'(x)$,
we have $\mbox{\boldmath $v$}=\mbox{\boldmath $A$}(x)\times\mbox{\boldmath $A$}'(x)\not=0$.
Since $\mbox{\boldmath $A$}(x)\cdot\mbox{\boldmath $v$}=\mbox{\boldmath $A$}'(x)\cdot\mbox{\boldmath $v$}=0$,
$\mbox{\boldmath $v$}$ satisfies
$$\displaystyle
S_x(\mbox{\boldmath $v$})
=\frac{\lambda(x)}{\sqrt{2(1+\mbox{\boldmath $A$}(x)\cdot\mbox{\boldmath $A$}'(x))}}\mbox{\boldmath $v$}.
$$
Next choose $\mbox{\boldmath $v$}'=\mbox{\boldmath $A$}(x)-\mbox{\boldmath $A$}'(x)$.
Since $\nu_x$ and $\mbox{\boldmath $A$}(x)+\mbox{\boldmath $A$}'(x)$ are parallel, $\mbox{\boldmath $v$}'\in T_xE_c(p,p')$.
We have
$$\displaystyle
S_x(\mbox{\boldmath $v$}')
=\frac{\lambda(x)\sqrt{1+\mbox{\boldmath $A$}(x)\cdot\mbox{\boldmath $A$}'(x)}}{2\sqrt{2}}\mbox{\boldmath $v$}'.
$$
Note also that $\mbox{\boldmath $v$}$ and $\mbox{\boldmath $v$}'$ are perpendicular to each other.
Therefore the eigenvalues of $S_x$ consists of two real numbers:
$$\begin{array}{c}
\displaystyle
k_1(x)=\frac{\lambda(x)}{\sqrt{2(1+\mbox{\boldmath $A$}(x)\cdot\mbox{\boldmath $A$}'(x))}},\\
\\
\displaystyle
k_2(x)=\frac{\lambda(x)\sqrt{1+\mbox{\boldmath $A$}(x)\cdot\mbox{\boldmath $A$}'(x)}}{2\sqrt{2}}.
\end{array}
\tag {A.3}
$$

If $\mbox{\boldmath $A$}(x)=\mbox{\boldmath $A$}'(x)$, then $\nu_x=\mbox{\boldmath $A$}(x)$.
Since $\mbox{\boldmath $v$}\cdot\nu_x=0$ for all $\mbox{\boldmath $v$}\in T_x(E_c(p,p')$, from (A.2)
we obtain
$$\displaystyle
S_x(\mbox{\boldmath $v$})
=\frac{\lambda(x)}{2}\mbox{\boldmath $v$},\,\,\forall\mbox{\boldmath $v$}\in\,T_x(E_c(p,p')).
$$
Thus the set of all eigenvalues of $S_x$ consists of only $\lambda(x)/2$.
Therefore $k_1(x)$, $k_2(x)$ given by (A.3) covers also this special
case.
Note that $k_2(x)\le k_1(x)$ and $k_1(x)=k_2(x)$ if and only
if $\mbox{\boldmath $A$}(x)=\mbox{\boldmath $A$}'(x)$. Therefore the Gauss curvature $K(x)$ at $x\in
E_c(p,p')$ and the mean curvature $H(x)$ with respect to $\nu_x$
are
$$\begin{array}{c}
\displaystyle
K(x)=k_1(x)k_2(x)=\frac{\lambda(x)^2}{4},\\
\\
\displaystyle
H(x)=\frac{k_1(x)+k_2(x)}{2}
=\frac{\lambda(x)}{8}(3+\mbox{\boldmath $A$}(x)\cdot\mbox{\boldmath $A$}'(x)).
\end{array}
$$

\subsection{Proof of Lemma 5.1}

Replacing $p'$ with $p'+s\mbox{\boldmath $A$}'$,
it follows from (4.18) and (4.19) that

$$\begin{array}{c}
\displaystyle
\text{det}\,
\left(S_q(E_{c-s}(p,p'+s\mbox{\boldmath $A$}'))-S_q(\partial D)\right)
=\text{det}\,\left(\frac{\lambda(q;p,p'+s\mbox{\boldmath $A$}')}{\sqrt{2(1+\mbox{\boldmath $A$}\cdot\mbox{\boldmath $A$}')}}
\delta_{kj}-M_{kj}\right)\\
\\
\displaystyle
=\left(\frac{\lambda(q;p,p'+s\mbox{\boldmath $A$}')}{\sqrt{2(1+\mbox{\boldmath $A$}\cdot\mbox{\boldmath $A$}')}}\right)^2
-\frac{\lambda(q;p,p'+s\mbox{\boldmath $A$}')}{\sqrt{2(1+\mbox{\boldmath $A$}\cdot\mbox{\boldmath $A$}')}}\text{Trace}\,(M_{kj})
+\text{det}\,(M_{kj}),
\end{array}
\tag {A.4}
$$
where
$$\displaystyle
M_{kj}
=\frac{\lambda(q;p,p'+s\mbox{\boldmath $A$}')}{2\sqrt{2(1+\mbox{\boldmath $A$}\cdot\mbox{\boldmath $A$}')}}
(\mbox{\boldmath $A$}\cdot\mbox{\boldmath $e$}_k\mbox{\boldmath $A$}\cdot\mbox{\boldmath $e$}_j
+\mbox{\boldmath $A$}'\cdot\mbox{\boldmath $e$}_k\mbox{\boldmath $A$}'\cdot\mbox{\boldmath $e$}_j)
+\frac{\partial^2 f}{\partial\sigma_k\sigma_j}(0).
\tag {A.5}
$$
From (4.14) we have
$$\displaystyle
\mbox{\boldmath $A$}\cdot\nu_q=\mbox{\boldmath $A$}'\cdot\nu_q
=-\sqrt{\frac{1+\mbox{\boldmath $A$}\cdot\mbox{\boldmath $A$}'}{2}}.
$$
These give
$$\displaystyle
\sum_{k=1}^2(
\vert\mbox{\boldmath $A$}\cdot\mbox{\boldmath $e$}_k\vert^2
+\vert\mbox{\boldmath $A$}'\cdot\mbox{\boldmath $e$}_k\vert^2)
=1-\mbox{\boldmath $A$}\cdot\mbox{\boldmath $A$}'.
$$
Thus one gets
$$\displaystyle
\text{Trace}\,(M_{kj})
=\frac{\lambda(q;p,p'+s\mbox{\boldmath $A$}')(1-\mbox{\boldmath $A$}\cdot\mbox{\boldmath $A$}')}{2\sqrt{2(1+\mbox{\boldmath $A$}\cdot\mbox{\boldmath $A$}')}}
+2H_{\partial D}(q).
\tag {A.6}
$$

For the computation of $\text{det}\,(M_{kj})$ we prepare the following formula.

\proclaim{\noindent Proposition A.2.}
Let $B$ be a $2\times 2$-matrix, $\mbox{\boldmath $c$}$ and $\mbox{\boldmath $c$}'$ be two-dimensional vectors.
Let $\gamma$ be a constant.
Let
$$\displaystyle
M=\gamma(\mbox{\boldmath $c$}\otimes\mbox{\boldmath $c$}+\mbox{\boldmath $c$}'\otimes\mbox{\boldmath $c$}')+B.
\tag {A.7}
$$
We have
$$\begin{array}{c}
\displaystyle
\text{det}\,M
\\
\\
\displaystyle
=\gamma^2\left(
\text{det}\,\left(\begin{array}{lr}
\displaystyle
c_1 & c'_1\\
\\
\displaystyle
c_2 & c'_2
\end{array}
\right)\right)^2
+\gamma\left(
B
\left(\begin{array}{c}
\displaystyle
c_2\\
\\
\displaystyle
-c_1
\end{array}
\right)
\cdot
\left(\begin{array}{c}
\displaystyle
c_2\\
\\
\displaystyle
-c_1
\end{array}
\right)
+
B
\left(\begin{array}{c}
\displaystyle
c'_2\\
\\
\displaystyle
-c'_1
\end{array}
\right)
\cdot
\left(\begin{array}{c}
\displaystyle
c'_2\\
\\
\displaystyle
-c'_1
\end{array}
\right)
\right)
+\text{det}\,B.
\end{array}
$$

\endproclaim

{\it\noindent Proof.}
We have
$$\begin{array}{c}
\displaystyle
\text{det}\,M\\
\\
\displaystyle
=(\gamma (c_1^2+c_1'^2)
+b_{11})(\gamma (c_2^2+c_2'^2)
+b_{22})
-(\gamma(c_1c_2+c_1'c_2')+b_{12})
(\gamma(c_2c_1+c_2'c_1')+b_{21})\\
\\
\displaystyle
=\gamma^2(c_1^2+c_1'^2)(c_2^2+c_2'^2)
+\gamma((c_1^2+c_1'^2)b_{22}+(c_2^2+c_2'^2)b_{11})
+b_{11}b_{22}\\
\\
\displaystyle
-\gamma^2(c_1c_2+c_1'c_2')(c_2c_1+c'_2c_1')
-\gamma((c_1c_2+c_1'c_2')b_{21}
+(c_2c_1+c_2'c_1')b_{12})-b_{12}b_{21}\\
\\
\displaystyle
=\gamma^2((c_1^2+c_1'^2)(c_2^2+c_2'^2)
-(c_1c_2+c_1'c_2')^2)\\
\\
\displaystyle
+\gamma((c_1^2+c_1'^2)b_{22}+(c_2^2+c_2'^2)b_{11}
-(c_1c_2+c_1'c_2')b_{21}
-(c_2c_1+c_2'c_1')b_{12})
+\text{det}\,B\\
\\
\displaystyle
=\gamma^2(c_1c_2'-c_1'c_2)^2
\\
\\
\displaystyle
+\gamma
((c_1^2+c_1'^2)b_{22}+(c_2^2+c_2'^2)b_{11}
-(c_1c_2+c_1'c_2')(b_{12}+b_{21}))+\text{det}\,B.
\end{array}
$$

\noindent
$\Box$

Note that $(M_{kj})$ given by (A.5) coincides with (A.7) in the
case when $\mbox{\boldmath $c$}=(\mbox{\boldmath
$A$}\cdot\mbox{\boldmath $e$}_j)$, $\mbox{\boldmath
$c$}'=(\mbox{\boldmath $A$}'\cdot\mbox{\boldmath $e$}_j)$,
$B=\nabla^2f(0)$ and
$$\displaystyle
\gamma=\frac{\lambda(q;p,p'+s\mbox{\boldmath $A$}')}{2\sqrt{2(1+\mbox{\boldmath $A$}\cdot\mbox{\boldmath $A$}')}}.
$$

From (4.15) we have
$$\displaystyle
\text{det}\,\left(\begin{array}{lr}
\displaystyle
\mbox{\boldmath $A$}\cdot\mbox{\boldmath $e$}_1 & \mbox{\boldmath $A$}'\cdot\mbox{\boldmath $e$}_1\\
\\
\displaystyle
\mbox{\boldmath $A$}\cdot\mbox{\boldmath $e$}_2 & \mbox{\boldmath $A$}'\cdot\mbox{\boldmath $e$}_2
\end{array}
\right)
=\text{det}\,\left(\begin{array}{lr}
\displaystyle
\mbox{\boldmath $A$}\cdot\mbox{\boldmath $e$}_1 & -\mbox{\boldmath $A$}\cdot\mbox{\boldmath $e$}_1\\
\\
\displaystyle
\mbox{\boldmath $A$}\cdot\mbox{\boldmath $e$}_2 & -\mbox{\boldmath $A$}\cdot\mbox{\boldmath $e$}_2
\end{array}
\right)=0;
\tag {A.8}
$$
from (4.19) we have
$$\begin{array}{c}
\displaystyle
\nabla^2f(0)
\left(\begin{array}{c}
\displaystyle
\mbox{\boldmath $A$}\cdot\mbox{\boldmath $e$}_2\\
\\
\displaystyle
-\mbox{\boldmath $A$}\cdot\mbox{\boldmath $e$}_1
\end{array}
\right)
\cdot
\left(\begin{array}{c}
\displaystyle
\mbox{\boldmath $A$}\cdot\mbox{\boldmath $e$}_2\\
\\
\displaystyle
-\mbox{\boldmath $A$}\cdot\mbox{\boldmath $e$}_1
\end{array}
\right)
\\
\\
\displaystyle
=S_{q}(\partial D)(
(\mbox{\boldmath $A$}\cdot\mbox{\boldmath $e$}_2)\mbox{\boldmath $e$}_1
-(\mbox{\boldmath $A$}\cdot\mbox{\boldmath $e$}_1)\mbox{\boldmath $e$}_2)
\cdot
((\mbox{\boldmath $A$}\cdot\mbox{\boldmath $e$}_2)\mbox{\boldmath $e$}_1
-(\mbox{\boldmath $A$}\cdot\mbox{\boldmath $e$}_1)\mbox{\boldmath $e$}_2).
\end{array}
\tag {A.9}
$$
By Lemma 4.3, we have
$$\displaystyle
\mbox{\boldmath $A$}\times\mbox{\boldmath $A$}'
=-\sqrt{2(1+\mbox{\boldmath $A$}\cdot\mbox{\boldmath $A$}')}
((\mbox{\boldmath $A$}\cdot\mbox{\boldmath $e$}_2)\mbox{\boldmath $e$}_1
-(\mbox{\boldmath $A$}\cdot\mbox{\boldmath $e$}_1)\mbox{\boldmath $e$}_2).
$$
This gives
$$\begin{array}{c}
\displaystyle
S_{q}(\partial D)(
(\mbox{\boldmath $A$}\cdot\mbox{\boldmath $e$}_2)\mbox{\boldmath $e$}_1
-(\mbox{\boldmath $A$}\cdot\mbox{\boldmath $e$}_1)\mbox{\boldmath $e$}_2)
\cdot
((\mbox{\boldmath $A$}\cdot\mbox{\boldmath $e$}_2)\mbox{\boldmath $e$}_1
-(\mbox{\boldmath $A$}\cdot\mbox{\boldmath $e$}_1)\mbox{\boldmath $e$}_2)
\\
\\
\displaystyle
=\frac{S_q(\partial D)(\mbox{\boldmath $A$}\times\mbox{\boldmath $A$}')
\cdot
(\mbox{\boldmath $A$}\times\mbox{\boldmath $A$}')}
{2(1+\mbox{\boldmath $A$}\cdot\mbox{\boldmath $A$}')}
\end{array}
$$
and thus from (A.9) one gets
$$\displaystyle
\nabla^2f(0)
\left(\begin{array}{c}
\displaystyle
\mbox{\boldmath $A$}\cdot\mbox{\boldmath $e$}_2\\
\\
\displaystyle
-\mbox{\boldmath $A$}\cdot\mbox{\boldmath $e$}_1
\end{array}
\right)
\cdot
\left(\begin{array}{c}
\displaystyle
\mbox{\boldmath $A$}\cdot\mbox{\boldmath $e$}_2\\
\\
\displaystyle
-\mbox{\boldmath $A$}\cdot\mbox{\boldmath $e$}_1
\end{array}
\right)
=
\frac{S_q(\partial D)(\mbox{\boldmath $A$}\times\mbox{\boldmath $A$}')
\cdot
(\mbox{\boldmath $A$}\times\mbox{\boldmath $A$}')}{2(1+\mbox{\boldmath $A$}\cdot\mbox{\boldmath $A$}')}.
\tag {A.10}
$$
Since
$$\displaystyle
(\mbox{\boldmath $A$}'\cdot\mbox{\boldmath $e$}_2)\mbox{\boldmath $e$}_1
-(\mbox{\boldmath $A$}'\cdot\mbox{\boldmath $e$}_1)\mbox{\boldmath $e$}_2
=-
((\mbox{\boldmath $A$}\cdot\mbox{\boldmath $e$}_2)\mbox{\boldmath $e$}_1
-(\mbox{\boldmath $A$}\cdot\mbox{\boldmath $e$}_1)\mbox{\boldmath $e$}_2),
$$
we obtain also
$$\displaystyle
\nabla^2f(0)
\left(\begin{array}{c}
\displaystyle
\mbox{\boldmath $A$}'\cdot\mbox{\boldmath $e$}_2\\
\\
\displaystyle
-\mbox{\boldmath $A$}'\cdot\mbox{\boldmath $e$}_1
\end{array}
\right)
\cdot
\left(\begin{array}{c}
\displaystyle
\mbox{\boldmath $A$}'\cdot\mbox{\boldmath $e$}_2\\
\\
\displaystyle
-\mbox{\boldmath $A$}'\cdot\mbox{\boldmath $e$}_1
\end{array}
\right)
=
\frac{S_q(\partial D)(\mbox{\boldmath $A$}\times\mbox{\boldmath $A$}')
\cdot
(\mbox{\boldmath $A$}\times\mbox{\boldmath $A$}')}{2(1+\mbox{\boldmath $A$}\cdot\mbox{\boldmath $A$}')}.
\tag {A.11}
$$
From Proposition A.2, (A.8), (A.10) and (A.11) we obtain
$$\displaystyle
\text{det}\,(M_{kj})
=\frac{\lambda(q;p,p'+s\mbox{\boldmath $A$}')}{\sqrt{2}
(1+\mbox{\boldmath $A$}\cdot\mbox{\boldmath $A$}')^{3/2}}S_q(\partial D)(\mbox{\boldmath $A$}\times\mbox{\boldmath $A$}')
\cdot
(\mbox{\boldmath $A$}\times\mbox{\boldmath $A$}')
+K_{\partial D}(q).
$$
Substituting this together with (A.6) into (A.4), we obtain (5.6).

\vskip1cm
\noindent
e-mail address

ikehata@math.sci.gunma-u.ac.jp

\end{document}